\documentclass[12pt]{amsart}
\usepackage{graphics}
\usepackage{amsthm}
\usepackage{amssymb}
\usepackage{amsmath}
\usepackage{amsfonts}
\usepackage{epic}
\usepackage{curves}
\usepackage{epsfig}
\usepackage{bm}
\usepackage{amscd}
\input{psfig.sty}

\input diagrams

\newcommand{\baseRing}[1]{\ensuremath{\mathbb{#1}}}
\newcommand{\Z}{\baseRing{Z}}
\newcommand{\C}{\baseRing{C}}
\newcommand{\N}{\baseRing{N}}
\newcommand{\R}{\baseRing{R}}
\newcommand{\Q}{\baseRing{Q}}
\newcommand{\PP}{\baseRing{P}}

\theoremstyle{plain}
\newtheorem{theorem}{Theorem}[section]
\newtheorem{lemma}[theorem]{Lemma}
\newtheorem{corollary}[theorem]{Corollary}
\newtheorem{proposition}[theorem]{Proposition}

\theoremstyle{definition}
\newtheorem{definition}[theorem]{Definition}

\newtheorem{remark}[theorem]{Remark}
\newtheorem{example}[theorem]{Example}

\def\supp{ {\rm supp}\, }

\def\rank{ {\rm rank}\, }
\def\vol{ {\rm vol}\, }

\def\ini{ {\rm in}\, }
\def\ind{ {\rm ind}\, }
\def\ker{ {\rm ker}\, }

\def\re{ {\rm Re}\, }

\def\Horn{ {\rm Horn}\, }

\title{Bivariate hypergeometric $D$-modules}

\author{Alicia Dickenstein}
\address{Dto. de Matem\'atica, FCEyN,
Universidad de Buenos Aires, \newline \indent
(1428) Buenos Aires, Argentina.}
\email{alidick@dm.uba.ar}
\thanks{Alicia Dickenstein was partially supported  by
UBACYT X052 and ANPCYT 03-6568,
Argentina.}

\author{Laura Felicia Matusevich}
\address{Department of Mathematics, Harvard University, \newline \indent
Cambridge, MA 02138, USA.}
\email{laura@math.harvard.edu}
\thanks{Laura Felicia Matusevich was partially supported by a Sarah~M. Hallam 
fellowship at UC~Berkeley, and Liftoff fellowship from the
Clay Mathematics Institute.}

\author{Timur Sadykov}
\address{Department of Mathematics, \newline
\indent The University of Western Ontario,  \newline
\indent London, Ontario N6A~5B7, Canada.}
\email{tsadykov@uwo.ca}
\thanks{Timur Sadykov was partially supported by the Russian Ministry of
Education, grant E~02-1-138.}

\begin{document}


\begin{abstract}
We undertake the study of bivariate Horn systems for generic
parameters. We prove that these hypergeometric systems are holonomic,
and we provide an explicit formula for their holonomic rank as well as
bases of their spaces of complex holomorphic solutions.
We also obtain analogous results for the generalized hypergeometric systems
arising from lattices of any rank.
\end{abstract}


\maketitle

\tableofcontents

\section{Introduction}
\label{sec:intro}

Classically, there have been two main directions in the study of
hypergeometric functions.
The first of these is to study the properties of a particular
series, analyze its convergence, compute its values at some
specific points providing combinatorial identities, give integral
representations, and find relations with other series of the same kind.
Here one could refer to well known works of Gauss and
Euler, for instance,~\cite{euler} and~\cite{gauss}.

The other classical avenue of research is to find a differential equation
that our hypergeometric function satisfies, and to study all the solutions
of that equation. This approach was pioneered by Kummer, who showed that the
Gauss hypergeometric function:
\begin{align*}
f(z) & =F[a,b;c;z] \\ & = 1 + \frac{ab}{c}\frac{z}{1!}+
\frac{a(a+1)b(b+1)}{c(c+1)}\frac{z^2}{2!}+
\frac{a(a+1)(a+2)b(b+1)(b+2)}{c(c+1)(c+2)}\frac{z^3}{3!}+\cdots
\end{align*}
satisfies the differential equation:
\[
z(1-z)\frac{d^2f}{d z^2}+(c-(1+a+b)z)\frac{d f}{d z}-abf=0.
\]
Kummer
went on to find all of the solutions of this equation (see~\cite{kummer}).
He constructed twenty-four (Gauss) series that,
whenever~$a$,~$b$ and~$c$ are not integers, provide
representations of two linearly independent solutions
to the Gauss equation, that are valid in any region
of the complex plane. Riemann also had a fundamental influence in this field
\cite{riemann}. For more historical details on hypergeometric functions, and a
comprehensive treatment of their classical theory, see \cite{slater}.

Both of these approaches have been tried for bivariate hypergeometric
series. In his article~\cite{erdelyi}, Erd\'elyi gives a complete set
of solutions for the following system of two hypergeometric
equations in two variables:
\begin{align*}
\big( x (\theta_x + \theta_y +a)(\theta_x+b)-\theta_x(\theta_x+\theta_y+c-1)
\big) f = & \; 0 \; , \\
\big( y (\theta_x + \theta_y +a)(\theta_y+b')-\theta_y(\theta_x+\theta_y+c-1)
\big) f = & \; 0 \; ,
\end{align*}
where~$\theta_x=x\frac{\partial}{\partial x}$ and
$\theta_y=y \frac{\partial}{\partial y}$.
This is the system of equations for Appell's function~$F_1$, and for generic
values of the parameters~$a$,~$b$,~$b'$ and~$c$,
Erd\'elyi constructs more than~$120$ fully supported series
solutions through contour integration. By a {\em fully supported series},
we mean a series such that the
convex hull of the exponents of the monomials appearing with nonzero
coefficient contains a full dimensional cone.
The holonomic rank of this system,
that is, the dimension of its space of complex holomorphic solutions
around a nonsingular point, is~$3$.

Another interesting system of two second order hypergeometric
equations in two variables is:
\begin{align*}
\big( x(2 \theta_x-\theta_y+a')(2\theta_x-\theta_y+a'+1)-
(-\theta_x+2\theta_y+a)\theta_x \big) f = & \; 0 \; , \\
\big( y(-\theta_x+2\theta_y+a)(-\theta_x+2\theta_y+a+1)-
(2\theta_x-\theta_y+a')\theta_y \big) f = & \; 0 \; .
\end{align*}
This is the system of equations for Horn's function~$G_3$, and its holonomic
rank is~$4$. Erd\'elyi notes that, in a neighborhood of a given point,
three linearly independent solutions of this system can be obtained through
contour integral methods. He also finds a fourth linearly independent
solution: the Puiseux monomial~$x^{-(a+2a')/3} y^{-(2a+a')/3}$. He remarks
that the existence of this elementary solution is puzzling, especially since it
cannot be expressed using contour integration, and offers no explanation
for its occurrence.

One of the goals of this article is to give a formula for the rank
of a system of two hypergeometric equations in two variables
when the parameters are generic (cf.~Theorem \ref{thm:main-theorem}). 
We will explain
why the system for Appell's~$F_1$ has rank~$3$ and why the very
similar system for Horn's~$G_3$ has rank~$4$.
We will also show that Puiseux polynomial solutions are
a commonplace phenomenon. Moreover, we will prove that these systems of hypergeometric
equations are holonomic for a generic choice of the parameters.

Our starting point are the ideas of Gel$'$fand, Graev and
Retakh~\cite{gelfand-graev-retakh} about the $\Gamma$-series associated with
lattices, and how
they relate to Horn series. Notice that $\Gamma$-series as defined
in~\cite{gelfand-graev-retakh} are fully supported, and they do not account for the
Puiseux polynomial solutions of Horn systems.

Holomorphic series solutions to a Horn system are equivalent to
solutions of corresponding hypergeometric recursions (see Section~\ref{sec:timur},
specifically equation  (\ref{horn-recurrencies})), thus our study of Puiseux polynomial 
solutions also characterizes the solutions to these recurrences that have finite support.

Finally, since we will be dealing
with lattices that are not necessarily saturated, we also need to
study the generalized hypergeometric systems associated with lattices
(more general than the $A$-hypergeometric systems of Gel$'$fand,
Kapranov and Zelevinsky). We show that, for generic parameters, these
systems are also holonomic, without restriction on the number
of variables or rank of the corresponding lattice, and prove the expected
formula for their generic holonomic rank.


\section{Multivariate hypergeometric systems}
\label{sec:formula}

In order to accommodate two different sets of variables,
we denote by~$D_n$ the Weyl algebra with generators $x_1,\dots, x_n$,
$\partial_{x_1},\dots,\partial_{x_n}$, and by~$D_m$ the Weyl algebra
whose generators are $y_1, \dots ,y_m$, $\partial_{y_1},\dots ,\partial_{y_m}$.
We set
$\theta_{x_j}=x_j \partial_{x_j}$ for $1 \leq j \leq n$, and
$\theta_{y_i}=y_i\partial_{y_i}$, for $ 1 \leq i \leq m$.
We also define $\theta_x=(\theta_{x_1}, \dots , \theta_{x_n})$ and
$\theta_y=(\theta_{y_1},\dots ,\theta_{y_m})$.
When the meaning is clear, we will drop many of the subindices to simplify the
notation.

We fix a matrix $A=(a_{ij}) \in \Z^{(n-m) \times n}$ of full rank~$n-m$
whose first row is the vector $(1,\dots ,1)$,
and a matrix ${\mathcal{B}} \in \Z^{n \times m}=(b_{ji})$ of full rank~$m$
such that $A \cdot {\mathcal{B}} = 0$. For $1 \leq j \leq m$, set
$b_j = (b_{j1},\dots , b_{jm}) \in \Z^m$ the~$j$-th row of~${\mathcal{B}}$.
The (positive) greatest common divisor of the maximal
minors of the matrix~${\mathcal{B}}$ is denoted by~$g$.

\smallskip

For $i = 1,\dots , m$, and a fixed parameter vector $c=(c_1,\dots, c_n) \in \C^n$,
we let
\begin{align}\label{def:horn}
&{\bm{P}}_i = \prod_{b_{ji}<0} \prod_{l=0}^{|b_{ji}|-1}( b_j \cdot \theta_y +
c_j-l), \\
&\bm{Q}_i = \prod_{b_{ji}>0} \prod_{l=0}^{b_{ji}-1} (b_j \cdot \theta_y
+ c_j-l) , \; \mbox{and} \\
& H_i = \bm{Q}_i-y_i\bm{P}_i ,
\end{align}
where $ b_j \cdot \theta_y = \sum_{k=1}^m b_{jk} \theta_{y_k}.$
The operators~$H_i$ are the {\em Horn operators} corresponding to the
lattice $L_{{\mathcal{B}}}=\{ {\mathcal{B}}\cdot z : z \in \Z^m \}$
and the parameter
vector~$c$. We call $d_i=\sum_{b_{ij}>0} b_{ij} = - \sum_{b_{ij}<0} b_{ij}$
the order of the operator~$H_i$.

\begin{definition}
The {\em Horn system} is the following left ideal of~$D_m$:
\[ \Horn({\mathcal{B}},c) = \langle H_1,\dots ,H_m \rangle \subseteq D_m .\]
\end{definition}
Now denote by~$b^{(i)}$ the columns of the matrix~$\mathcal{B}$.
Any vector~$u \in \R^n$ can be written as $u=u_+-u_-$, where $(u_+)_i = \max(u_i,0)$,
and $(u_-)_i= - \min(u_i,0)$.
For $i=1,\dots ,m$, we let:
\[
T_i = \partial_x^{b^{(i)}_+} - \partial_x^{b^{(i)}_-},
\]
here we use multi-index notation $\partial_x^v=\partial_{x_1}^{v_1}\cdots
\partial_{x_n}^{v_n}$.
More generally, for any $u \in L_{{\mathcal{B}}}$, set
\[ T_u = \partial_x^{u_+}-\partial_x^{u_-} . \]
These are the {\em lattice operators} arising from~$L_{{\mathcal{B}}}$.

\begin{definition}
The {\em lattice ideal} arising from~$L_{{\mathcal{B}}}$ is:
\[I_{{\mathcal{B}}}
= \langle T_u : u \in L_{{\mathcal{B}}} \rangle \subseteq \C[\partial_{x_1},
\dots ,\partial_{x_n} ] .\]
Recall that the {\em toric ideal} corresponding to~$A$ is:
\[ I_A  = \langle T_u : u \in \ker_{\Z}(A) \rangle
\subseteq \C[\partial_{x_1}, \dots ,\partial_{x_n} ] .\]
We will also denote:
\[ I = \langle T_1,\dots ,T_m \rangle
\subseteq \C[\partial_{x_1}, \dots ,\partial_{x_n} ] .\]
The ideal~$I$ is called a {\em lattice basis ideal}.
Notice that, for~$m=2$,~$I$ is a complete intersection. This is not
necessarily true if~$m>2$.
\end{definition}
Lattice ideals and toric ideals have been extensively studied (see,
for instance,~\cite{binomialideals},~\cite{gb&cp}).
Lattice basis ideals were introduced in~\cite{latticebasis}.

There is a natural system of differential equations arising from a toric 
ideal~$I_A$ and a parameter vector. This system, called the
{\em $A$-hypergeometric system with
parameter~$A\cdot c$}, is defined as:
\[ H_A(A\cdot c) = I_A + \langle \sum_{j=1}^n a_{ij}x_j \partial_{x_j} -
(A\cdot c)_i : i=1,\dots ,n-m \rangle \subseteq D_n .\]
{}From now on we will use the notation $\langle A\cdot \theta - A\cdot c \rangle$
to mean $\langle \sum_{j=1}^n a_{ij}x_j \partial_{x_j} -
(A\cdot c)_i : i=1,\dots ,n-m \rangle$.

$A$-hypergeometric systems were first defined by Gel$'$fand, Graev and Zelevinsky
in~\cite{GGZ}, and their systematic analysis was started by Gel$'$fand,
Kapranov and Zelevinsky (see, for instance,~\cite{GKZ}). Saito, Sturmfels and
Takayama have used Gr\"{o}bner deformations in the Weyl algebra
to study $A$-hypergeometric systems
(see~\cite{SST}). In this article, we will extend this approach to the case
of Horn systems.

Gel$'$fand, Graev and Retakh have also considered the {\em hypergeometric
system associated with the lattice
$L_{{\mathcal{B}}} = \{ {\mathcal{B}} \cdot z : z \in \Z^m \}$}, which is
defined to be the left $D_n$-ideal:
\[ I_{{\mathcal{B}}} + \langle A\cdot \theta - A\cdot c \rangle \subseteq D_n.\]
We now introduce the $D_n$-ideal~$H_{{\mathcal{B}}}(c)$, that
is very closely related to the Horn system $\Horn({{\mathcal{B}}},c)$:
\[H_{{\mathcal{B}}}(c) = I +
\langle A \cdot \theta - A\cdot c \rangle \subseteq D_n .\]
The results in Section~\ref{sec:translation} imply that, for 
generic~$c$, there is a vector space
isomorphism between the solution spaces of $\Horn({\mathcal{B}},c)$
and~$H_{{\mathcal{B}}}(c)$. Thus, we have two points of view to study
Horn hypergeometric functions. We also call~$H_{{\mathcal{B}}}(c)$
a Horn system, when the context is clear.

\begin{remark}
We have defined the Horn operators using falling factorials because
this formulation will make clearer the relationship between
$\Horn({\mathcal{B}},c)$ and~$H_{{\mathcal{B}}}(c)$, but it is just as
legal to define Horn systems using rising factorials, as it is done
in many classical sources. For instance, the Horn
and Appell systems from the previous section naturally lend themselves
to a rising factorial formulation. This is not really a difficulty, since
switching between rising and falling factorials in the
definition of Horn systems is a matter of shifting the parameters by integers.
\end{remark}

It is a well known result of Adolphson~\cite{adolphson} that, for
generic parameters~$A\cdot c$, the holonomic rank of the $A$-hypergeometric
system equals the normalized volume~$\vol(A)$ of the convex hull of the
columns of~$A$, which is also the degree of the toric ideal~$I_A$.
Our goal is to obtain an explicit expression in this spirit for bivariate
Horn systems. Previous work in this direction required very strong
assumptions (see~\cite{timur}).

\begin{definition}\label{def:index} In the case that~$m=2$, we set
\[ \nu_{ij} = \left\{ \begin{array}{ll} \min(|b_{i1}b_{j2}|,|b_{j1}b_{i2}|),
\quad  & \mbox{if~$b_i$,~$b_j$ are in the interior of opposite quadrants of~$\Z^2$}, \\
0 & \mbox{otherwise}, \end{array} \right.\]
for $1 \leq i,j\leq n$. The number~$\nu_{ij}$ is called the {\em index} associated
to~$b_i$ and~$b_j$.
\end{definition}

The following is the main result in this article, which follows from 
Corollary~\ref{corol:formula} and 
Theorems~\ref{thm:hbc-is-holonomic},~\ref{thm:rank-formula}, and~\ref{hornholonom}.

\begin{theorem}
\label{thm:main-theorem}
Let
${\mathcal{B}}$ be an~$n \times 2$ integer matrix of full rank such that
its rows $b_1,\dots, b_n$ satisfy $b_1+\cdots+b_n=0$.
If $c\in \C^n$ is a generic parameter vector, then the ideals
$\Horn({\mathcal{B}},c)$ and $H_{{\mathcal{B}}}(c)$ are holonomic. Moreover,
\[ \rank(H_{{\mathcal{B}}}(c)) = \rank(\Horn({\mathcal{B}},c)) = d_1d_2 - 
\sum_{\begin{array}{c} b_i,b_j \\ \mbox{\tiny dependent} \end{array}} \nu_{ij} \;
 = \, g \cdot \vol(A)+
 \sum_{\begin{array}{c} b_i,b_j \\  \mbox{\tiny independent} \end{array}} \nu_{ij} \; ,\]
where the first summation runs over linearly dependent pairs~$b_i$,~$b_j$
of rows of~${\mathcal{B}}$
that lie in opposite open quadrants of~$\Z^2$, and the second summation
runs over linearly independent such pairs.
\end{theorem}

We can also give an explicit basis for the solution space
of $\Horn({{\mathcal{B}}},c)$ (and of~$H_{{\mathcal{B}}}(c)$)
(Theorem~\ref{thm:sol-basis}), and compute the
exact dimension of the subspace of Puiseux polynomial solutions (Theorem
\ref{thm:puiseux}).


\section{Some observations about Horn systems}

The Horn system $\Horn({\mathcal{B}},c)$
is always compatible, even if~$c$ is not generic, in the sense that
its solution space is always nonempty. First of all, the constant zero function
is always a solution of $\Horn({\mathcal{B}},c)$, since this system is homogeneous.
Moreover, as we will see in Section~\ref{sec:translation}, all the solutions
of the $A$-hypergeometric system~$H_A(A\cdot c)$ are solutions of
$H_{{\mathcal{B}}}(c)$, and these can be transformed into solutions of
$\Horn({\mathcal{B}},c)$ (see
Corollary~\ref{coro:vector-space-iso}), so that, under the assumptions that
${{\mathcal{B}}}$ is~$n\times m$ of full rank~$m$,~$n>m$, with all column sums equal
to zero, $\Horn({\mathcal{B}},c)$ always has nonzero solutions, since
$H_A(A\cdot c)$ always has nonzero solutions (its solution space has
dimension at least $\deg(I_A)=\vol(A)$, see~\cite[Theorem 3.5.1]{SST}).

It is easy to understand how the Horn system $\Horn({\mathcal{B}},c)$
changes
if we choose a new parameter vector~$c'$, as long as $A\cdot c' = A\cdot c$.
As a matter of fact, if $c = c' + {\mathcal{B}}\cdot z$, for some
$z \in \C^m$, then it is easy to see that~$f(y)$ is a solution of
$\Horn({\mathcal{B}},c')$ if and only if~$y^z f(y)$ is a solution of
$\Horn({\mathcal{B}},c)$. Notice also that the system~$H_{{\mathcal{B}}}(c)$
depends only on~$A\cdot c$, so that $H_{{\mathcal{B}}}(c)=H_{{\mathcal{B}}}(c')$
if $A\cdot c = A\cdot c'$.

A change in~$A\cdot c\ $ can, instead, alter the solution space of
$\Horn({\mathcal{B}},c)$ (and~$H_{{\mathcal{B}}}(c)$)
in dramatic ways. For instance, it could become
infinite-dimensional, as the following example shows.

\begin{example}

The Horn system defined by the operators
\begin{equation}
(\theta_{y_{1}} + \theta_{y_{2}} + c_{1}) \theta_{y_{i}} -
y_{i} (\theta_{y_{1}} + \theta_{y_{2}} + c_{2})
(\theta_{y_{1}} + \theta_{y_{2}} + c_{3}), \,\,\,\,\, i=1,2
\label{eqn:nonholo2dim}
\end{equation}
is not holonomic if $(c_{1} - c_{2})(c_{1} - c_{3})=0.$
Indeed, a holonomic system of equations can only have
a finite-dimensional space of analytic solutions.
However, since for $(c_{1} - c_{2})(c_{1} - c_{3})=0$
the operator $\theta_{y_{1}} + \theta_{y_{2}} + c_{1}$ can be factored out of
each of the operators in~(\ref{eqn:nonholo2dim}),
it follows that any function which is
annihilated by $\theta_{y_{1}} + \theta_{y_{2}} + c_{1}$ is a solution
to~(\ref{eqn:nonholo2dim}). Thus for any smooth univariate function~$u$
the product~$y_{2}^{-c_{1}}u(y_{1}/y_{2})$
satisfies~(\ref{eqn:nonholo2dim}).

Notice that for generic values of the parameters $c_{1},c_{2},c_{3}$
the system~(\ref{eqn:nonholo2dim}) is holonomic. One of its solutions is given
by the Gauss function $F[c_{2},c_{3};c_{1}; y_{1} + y_{2}].$ Of course, similar examples
can be given in any dimension.
\end{example}

We could also ask what happens if we choose another matrix~${\mathcal{B}}'$
such
that $A \cdot {\mathcal{B}}' = 0$.
Even if~$g=g'=1$, so that~${\mathcal{B}}$ and~${\mathcal{B}}'$ are two
Gale duals of~$A$, the associated Horn systems could have different
holonomic rank, as we see in Example~\ref{ex:diff-ranks}.
The systematic analysis of this question, in the case when~$m=2$ is one of
the main objectives of this article.

\begin{example}

\label{ex:diff-ranks}
We choose:
\[ A = \left( \begin{array}{cccc} 1 & 1 & 1 & 1 \\ 0 & 1 & 2 & 3
\end{array} \right) \; , \;
B = \left( \begin{array}{rr} 1 & 0 \\ -2 & 1 \\ 1 & -2 \\ 0 & 1 \end{array}
\right) \; , \;
B' = \left( \begin{array}{rr} 1 & 2 \\ -2 & -3 \\ 1 & 0 \\ 0 & 1 \end{array}
\right) .\]
Then, if~$c$ is a generic parameter vector,
$\rank(\Horn({\mathcal{B}},c))=4$, and
$\rank(\Horn({\mathcal{B}}',c))=6$, as a consequence of Theorem~\ref{thm:main-theorem}.
This can be verified for specific values of~$c$ using the
computer algebra system {\sl Macaulay 2}~\cite{M2}. However, by Theorem
\ref{thm:splitting-series},  these
two hypergeometric systems share all fully supported solutions.
\end{example}

Notice that the definition of $\Horn({{\mathcal{B}}},c)$ makes sense even
if~${{\mathcal{B}}}$ is a square matrix, or if the rows of~${{\mathcal{B}}}$ do not
add up to zero, or even if~${{\mathcal{B}}}$ does not have full rank.
As a matter of fact, we will need to consider such Horn systems on our way to proving
results about the case when~${{\mathcal{B}}}$ is~$n\times m$ of full rank~$m$, 
$m<n$, and the rows of~${{\mathcal{B}}}$ add up to zero.
Many of the examples will also concern Horn systems with~$n=m$. We remark that
if~${{\mathcal{B}}}$ is square and nonsingular, then~$H_{{\mathcal{B}}}(c)$ is a system
of differential equations with constant coefficients, not depending on~$c$.


\section{Preliminaries on codimension~$2$ binomial ideals}

\label{sec:binomialideals}

In this section we collect some results about lattice ideals and
lattice basis ideals that will be necessary to study Horn systems.
Although this section is about commutative algebra, our indeterminates
will be called $\partial_1,\dots,\partial_n$ for consistency with the
notation for differential equations.

Recall that ${{\mathcal{B}}}=(b_{ji})$ is an~$n\times m$ integer matrix of
full rank~$m$ with all column sums equal to zero. The following ideal
is called a {\em lattice ideal}:
\[ I_{{\mathcal{B}}} = \langle \partial^{u_+}-\partial^{u_-} : u=u_+-u_- \in
L_{{\mathcal{B}}} \rangle \subset \C[\partial_1,\dots ,\partial_n], \]
where $L_{{\mathcal{B}}} = \{ {\mathcal{B}}\cdot z : z \in \Z^m\}$ is the
rank-$m$ lattice spanned by the columns of~${{\mathcal{B}}}$.
For the purpose of this section, we could use any field of characteristic~$0$ 
instead of~$\C$, but later on, when we talk about complex holomorphic
solutions of differential equations, we will need our field to be the complex numbers.
We let~$A$ be any $(n-m)\times n$ integer matrix such that $A\cdot {\mathcal{B}}=0$. Then
the saturation of~$L_{{\mathcal{B}}}$ is the lattice $L=\ker_{\Z}(A)$.
Notice that the order of the group~$L/L_{{\mathcal{B}}}$ is~$g$,
the positive greatest common divisor of the maximal minors
of~${{\mathcal{B}}}$.

The ideal~$I_{{\mathcal{B}}}$ is homogeneous with respect to the usual
$\Z$-grading and hence defines a subscheme~$X_{{\mathcal{B}}}$ of~$\PP^{n-1}$.
Moreover, the ideal~$I_{{\mathcal{B}}}$ is always radical
and~$X_{{\mathcal{B}}}$ is the equidimensional union of $g=|L/L_{{\mathcal{B}}}|$
torus translates of the toric variety~$X_A$ defined by the reduced scheme
associated to~$L$ as above. This is deduced from~\cite{binomialideals} since
$(I_{{\mathcal{B}}}:\langle \partial_1,\dots ,\partial_n \rangle^{\infty})=I_{{\mathcal{B}}}$,
that is, no component of~$X_{{\mathcal{B}}}$ is contained in a
coordinate hyperplane.

These  torus translates can be described in terms of the order~$g$ group~$G_{{\mathcal{B}}}$
of all partial characters $\rho : L \rightarrow \C^*$ which extend the trivial
character $1:L_{{\mathcal{B}}}\rightarrow \C^*$, i.e.,~$\rho$ satisfying
$\rho(\ell+\ell')=\rho(\ell)\rho(\ell'), \forall \ \ell,\ell' \in L$
and $\rho(\ell)=1$, $\forall \ \ell \in L_{{\mathcal{B}}}$.

\begin{example}

We illustrate the previous decomposition in an example before writing
it down in general. Let
\[ \mathcal{{B}} = \left( \begin{array}{rr} -1 & 2 \\ 0 & -3 \\ 3 & 0 \\ -2 & 1
\end{array} \right) , \quad A=\left( \begin{array}{rrrr} 1 & 1 & 1 & 1 \\
0 & 1 & 2 & 3 \end{array} \right). \]
In this case~$g=3$. The scheme~$X_A$ is the twisted cubic, that is,
the closure of the torus orbit of the point $p_0 = (1:1:1:1) \in \PP^3$
under the torus action:
\begin{equation}
\label{eqn:torus-action}
\lambda \cdot (\partial_1:\partial_2:\partial_3:\partial_4) =
(\lambda^0\partial_1:\lambda^1\partial_2:\lambda^2\partial_3:\lambda^3\partial_4),
\quad \lambda \in \C^* .
\end{equation}

The group~$G_{{\mathcal{B}}}$ has order~$3$ and is isomorphic to the group
of cubic roots of unity $\{1,\omega,\omega^2\}$, where
$\omega=e^{\frac{2\pi i}{3}}$. Set $p_1=(1:1:\omega:1)$,
$p_2=(1:1:\omega^2:1)$ and denote by~$X_0$,~$X_1$ and~$X_2$ the respective
closure of the torus orbit under the action (\ref{eqn:torus-action}) of
$p_0$,~$p_1$ and~$p_2$. In particular,~$X_0=X_A$. Then
\[X_B = X_0 \cup X_1 \cup X_2 \]
and~$X_i$ is the image of~$X_0$ under the coordinatewise multiplication by~$p_i$,
$i=1,2$. Note that
\[ X_i = \{ (\partial_1:\dots :\partial_4):
\partial_1\partial_3-\omega^i \partial_2^2=\partial_3^2-\omega^{2i}\partial_2\partial_4
= \partial_2\partial_3-\omega^i \partial_1\partial_4 = 0 \} \]
so that the equations defining~$X_i$ are ``translations'' of the equations
for~$X_0=X_A$.
\end{example}

This can be phrased in general as follows:
Given $\rho \in G_{{\mathcal{B}}}$, let~$X_{\rho}$ denote zero scheme of the ideal:
\[ I_{\rho} = \langle \partial^{u_+}-\rho(u)\partial^{u_-}: u = u_+-u_- \in L\rangle .\]
Then the ideals~$I_{\rho}$ are prime, their intersection gives~$I_{{\mathcal{B}}}$ and
$X_{{\mathcal{B}}}=\cup_{\rho \in G_{{\mathcal{B}}}} X_{\rho}$. We refer
to~\cite{binomialideals} for a proof of these facts.

\medskip

Consider now the case~$m=2$ and recall that the lattice basis ideal
associated to~${\mathcal{B}}$ is the ideal
\[ I = \langle \partial^{u_+}-\partial^{u_-} : u \; \mbox{is a column of}\; {{\mathcal{B}}}
\rangle .\]
Its zero set consists of the union of~$X_{{\mathcal{B}}}$ with components that lie inside
coordinate hyperplanes. The following proposition, whose proof can be
found in~\cite{elim-codim-2}, gives the precise
primary decomposition of the ideal~$I$.
Denote $b_1,\dots,b_n \in \Z^2$ the row vectors of~${{\mathcal{B}}}$.
Let~$\nu_{ij}$ be the index associated to~$b_i$ and~$b_j$ as in 
Definition~\ref{def:index}.

\begin{proposition}

\label{propo:lattice-basis-description}
The ideal~$I$ has the following primary decomposition:
\[I = \big(\cap_{\rho \in G_{{\mathcal{B}}}} I_{\rho} \big)\cap
\big(\cap_{\nu_{ij}>0} I_{ij} \big)\]
where $\sqrt{I_{ij}} = \langle \partial_i, \partial_j \rangle$, and
the multiplicity of each~$I_{ij}$ is~$\nu_{ij}$, in the sense that
\[\dim_K (\C[\partial_1,\dots,\partial_n]/I_{ij})_{\langle
\partial_1,\dots,\hat{\partial_i},\dots,\hat{\partial_j},
\dots,\partial_n \rangle}=\nu_{ij},\]
where $K=\C(\partial_1,\dots,\hat{\partial_i},\dots,\hat{\partial_j},
\dots,\partial_n)$.
\end{proposition}

We then have

\begin{corollary}
\label{corol:formula}
For~$d_1$,~$d_2$ the degrees of the generators of~$I$,
\begin{equation}
\label{coro:formula}
d_1 \cdot d_2 - \sum_{b_i,b_j \; \mbox{\tiny
 dependent}} \nu_{ij} = g \cdot \vol(A)+
 \sum_{b_i,b_j \; \mbox{\tiny independent}} \nu_{ij} \; ,
\end{equation}
where the first summation runs over linearly dependent pairs~$b_i$,~$b_j$ 
of rows of~${\mathcal{B}}$
that lie in opposite open quadrants of~$\Z^2$, and the second summation
runs over linearly independent such pairs.
\end{corollary}

\begin{proof}
The degree of the complete intersection~$I$ is~$d_1 d_2$. By Proposition
\ref{propo:lattice-basis-description}, this number equals
\[ g \cdot \deg(I_A) + \sum \nu_{ij} \;, \]
where the sum runs over all pairs of rows of~${{\mathcal{B}}}$ in
opposite open quadrants of~$\Z^2$. Now the result follows from the fact that the
degree of~$I_A$ is exactly the normalized volume~$\vol(A)$ of the polytope
obtained by taking the convex hull of the columns of~$A$~\cite[Theorem 4.16]{gb&cp}.
\end{proof}

The following is another result related to the primary decomposition
of~$I$.

\begin{proposition}

\label{propo:alicia's-lemma}
Let~${{\mathcal{B}}} \in \Z^{n \times 2}$ of rank~$2$, with
rows $b_1,\dots ,b_n$, that add up to zero, and
$I_{{\mathcal{B}}}$,~$I$, the lattice and lattice basis ideals
associated to~${{\mathcal{B}}}$. For each $1 \leq i,j \leq n$,
$\nu_{ij}$ is as in Definition~\ref{def:index}. Set
\[ \alpha_i = \left\{ \begin{array}{ll} \max_j \; \nu_{ij} \;&
  \mbox{if} \; b_{i1}>0, \\
0 & \mbox{otherwise.}

\end{array} \right. \]
Then
\[ \partial^{\alpha} I_{{\mathcal{B}}} \subseteq I. \]
\end{proposition}

\begin{proof}

By Proposition~\ref{propo:lattice-basis-description}, it is enough to
prove that $\partial^{\alpha} \in \cap_{\nu_{ij}>0} I_{ij}$.
Assume that~$\nu_{ij}>0$. Then~$b_i$ and~$b_j$ lie in the interior of
opposite quadrants, so that either~$b_{i1}$ or~$b_{j1}$ is
positive, say $b_{i1}>0$, so that $\alpha_i \geq \nu_{ij}$.
We will be done if we show that $ \partial_{i}^{\nu_{ij}} \in I_{ij} .$
To do this, let~$\tilde{I}_{ij}$ be the localization of~$I_{ij}$ at $\langle
\partial_1,\dots,\hat{\partial_i},\dots,\hat{\partial_j},
\dots,\partial_n \rangle$ so that~$\tilde{I}_{ij}$ is an
artinian ideal of multiplicity~$\nu_{ij}$ in $K[\partial_i,\partial_j]$,
where $K=\C(\partial_1,\dots,\hat{\partial_i},\dots,\hat{\partial_j},
\dots,\partial_n)$.
Notice that, since $\# \{1,\partial_i,\dots
,\partial_i^{\nu_{ij}} \} = \nu_{ij} + 1$, these monomials must be
linearly dependent modulo~$\tilde{I}_{ij}$, so we can find
$g_0,\dots,g_{\nu_{ij}} \in K$ such that
\[ g_0 + g_1 \partial_i + \cdots + g_{\nu_{ij}} \partial_i^{\nu_{ij}}
\in \tilde{I}_{ij} .\]
But the radical of~$\tilde{I}_{ij}$ is $\langle \partial_i, \partial_j
\rangle$, so that $g_0 = 0$. Let $l = \min_{1\leq k \leq \nu_{ij}} \{
g_k \neq 0 \}$. Then, clearing denominators, we can find polynomials
$f_l,\dots , f_{\nu_{ij}}$ {\em not} involving the variables
$\partial_i, \partial_j$, $f_l \neq 0$, such that
\[ \partial_i^l(f_l +\cdots + f_{\nu_{ij}}\partial_i^{\nu_{ij}-l}) \in
I_{ij} .\]
Now, since~$I_{ij}$ is primary to $\langle \partial_i,\partial_j
\rangle$, and no power of $f_l +\cdots +
f_{\nu_{ij}}\partial_i^{\nu_{ij}-l}$ belongs to $\langle \partial_i,\partial_j
\rangle$, then~$\partial_i^l$ must belong to~$I_{ij}$. Since $l \leq
\nu_{ij}$, we are done.
\end{proof}

It is an interesting fact that the multiplicities of some of the components 
of~$I$ do not go down under Gr\"{o}bner deformation. Given $w \in \Z^n$,
and $f = \sum f_{\alpha}x^{\alpha}$ a homogeneous polynomial in
$\C[\partial_1,\dots ,\partial_n]$, let
\[ \ini_w(f) = \sum_{w\cdot \alpha \; \mbox{\tiny maximal over}\; f_{\alpha}\neq 0}
f_{\alpha}x^{\alpha}\]
and define
\[\ini_w(I) = \langle \ini_w(f): f \in I \setminus \{0\}\rangle. \]
The ideal~$\ini_w(I)$ is called the {\em initial ideal of~$I$ with respect to
the weight vector~$w$}. It is a monomial ideal
if~$w$ is generic (see~\cite{IVA} and~\cite[Chapter 15]{commalg}
for more on initial ideals, especially how to compute them).

\begin{lemma}

\label{lemma:multiplicities}
Let~$b_k$ and~$b_l$ be two linearly dependent rows of~${\mathcal{B}}$ lying
in opposite open quadrants of~$\Z^2$. If~$w$ is a generic weight vector, then
the multiplicity of the ideal
$\langle \partial_k , \partial_l \rangle$ as an associated prime of
$\ini_w(I)$ is the index~$\nu_{kl}$.
\end{lemma}

This proof was suggested to us by Ezra Miller, to whom we are very grateful.

\begin{proof}

Recall that the initial variety of~$\mathcal{V}(I)$
is the flat limit of a family that is obtained by a one parameter subgroup of the
torus acting on the zero set~$\mathcal{V}(I)$. 
The monomial components of~$\mathcal{V}(I)$
are invariant under this action, so in the limit, the only way that the multiplicity of
$\langle \partial_l, \partial_k \rangle$ could go up is if this prime is associated
to $\ini_w(I_{{\mathcal{B}}})$. Now, if~$b_k$ and~$b_l$ are linearly dependent,
$\langle \partial_k, \partial_l \rangle$ is not associated to
$\ini_w(I_{{\mathcal{B}}})$, this follows from the same arguments that proved
\cite[Lemma 2.3]{shared}.

\end{proof}


\section{$A$-hypergeometric solutions of the Horn system}

\label{sec:translation}

In this section we study the solutions of the
Horn system $\Horn({\mathcal{B}},c)$ that arise
from the $A$-hypergeometric system~$H_A(A\cdot c)$. Here, we do not
use the assumption that~$m=2$.
Recall that ${\mathcal{B}}=(b_{ji})$ is an rank~$m$ integer $n\times m$ matrix whose rows
add up to zero, and whose columns are denoted $b^{(1)},\dots ,b^{(m)}$ and let
$A=(a_{ij})$ be any rank $(n-m)$ integer $(n-m)\times n$ matrix
such that $A \cdot {\mathcal{B}}=0$. Here we assume that~$n>m$.

Consider the surjective map
\begin{align*}
x^{\mathcal{B}} :  (\C^*)^n &\rightarrow (\C^*)^m , \\
          x &\mapsto (\prod_{j=1}^n x_j^{b_{j1}},\dots ,
\prod_{j=1}^n x_j^{b_{jm}}) = (x^{b^{(1)}},\dots , x^{b^{(m)}}) .
\end{align*}

This map is open in the sense that it takes open sets to open sets.
We use it to relate the operators~$T_i$ in~$n$ variables
and the operators~$H_i$ in~$m$ variables, defined in Section~\ref{sec:formula}.

\begin{lemma}

\label{lemma:how-to-translate}
Let $U \subseteq (\C^*)^n$ be a simply connected open set and let
$V = x^{{\mathcal{B}}}(U)$.
We choose~$U$ small enough so that~$V$ is also simply connected.
Given a holomorphic function $\psi \in {\mathcal{O}}(V)$,
call $\varphi = x^c \psi(x^{{\mathcal{B}}})$. Then
\begin{enumerate}

\item
\label{it-1}

$\big(\sum_{j=1}^n a_{kj}x_j \partial_{x_j} \big)(\varphi) =
(A \cdot c)_k \varphi$, for $k=1, \dots ,n-m$.
\item
\label{it-2}

$T_i (\varphi) = 0 $ for $i=1,\dots ,m$ if and only if
$H_i(\psi) = 0$ for $i=1, \dots ,m$.
\item
\label{it-3}

Moreover, for any $u={\mathcal{B}}\cdot z  \in L_{{\mathcal{B}}}$, and
\[ H_u = \prod_{u_j>0} \prod_{l=0}^{u_j-1} (b_j \cdot \theta_y + c_j -l)
- y^z \prod_{u_j < 0} \prod_{l=0}^{|u_j|-1}(b_j\cdot \theta_y +c_j-l), \]
we have $T_u(\varphi)=0$ if and only if~$H_u(\psi)=0$.
\end{enumerate}

\end{lemma}

\begin{proof}

The verifications of the three assertions are very similar.
The main ingredients are the following identities:
\begin{equation}
\label{eqn:use-1}
\theta_{x_i} x^c = x^c ( \theta_{x_i}+c_i), \quad \mbox{(in~$D_n$)} ,
\end{equation}

\begin{equation}
\label{eqn:use-2}
\theta_{x_i}(\psi(x^{{\mathcal{B}}}))(x) =
\big[ (b_i \cdot \theta_y) \psi \big] (x^{{\mathcal{B}}}) ,
\end{equation}
which are easily checked.
Let us prove (\ref{it-2}). Call
$\tilde{T}_i = \prod_{b_{ji}>0} x_j^{b_{ji}} T_i$.
We have:
\begin{equation}
\label{eqn:star}
\tilde{T}_i = \prod_{b_{ji}>0} x_j^{b_{ji}} \prod_{b_{ji}>0} \partial_{x_j}^{b_{ji}}-
(x^{{\mathcal{B}}})_i \prod_{b_{ji}<0} x_j^{-b_{ji}} \prod_{b_{ji}<0}
\partial_{x_j}^{-b_{ji}} .
\end{equation}

Recall that $(x^{{\mathcal{B}}})_i = \prod_{j=1}^n x_j^{b_{ji}}$.
Using the identity:
\[ x^{\alpha} \partial_x^{\alpha} = \prod_{j=1}^n \prod_{l=0}^{\alpha_j-1}
(\theta_{x_j}-l)\; ,\]
equation (\ref{eqn:star}) is transformed into:
\[ \tilde{T}_i = \prod_{b_{ji}>0} \prod_{l=0}^{b_{ji}-1}(\theta_{x_j}-l) -
(x^{{\mathcal{B}}})_i \prod_{b_{ji}<0} \prod_{l=0}^{-b_{ji}-1} (\theta_{x_j}-l)
\]
Using (\ref{eqn:use-1}),
\begin{align*}
\tilde{T}_i(\varphi)& = \tilde{T}_i(x^c \psi(x^{{\mathcal{B}}})) \\
                      & = x^c\bigg(
\prod_{l=0}^{b_{ji}-1}(\theta_{x_j}+c_j-l) -
(x^{{\mathcal{B}}})_i \prod_{b_{ji}<0} \prod_{l=0}^{-b_{ji}-1}
(\theta_{x_j}+c_j-l) \bigg) (\psi(x^{{\mathcal{B}}})).
\end{align*}

Now (\ref{eqn:use-2}) implies that
\begin{align*}
\tilde{T}_i(\varphi) & = x^c\bigg(
\prod_{l=0}^{b_{ji}-1}(b_j \cdot \theta_{y}(\psi)+c_j-l) - \\
& \qquad - (x^{{\mathcal{B}}})_i \prod_{b_{ji}<0} \prod_{l=0}^{-b_{ji}-1}
(b_j \cdot \theta_{y}(\psi)+c_j-l) \bigg)(x^{{\mathcal{B}}}) \\
& = x^c H_i (\psi)((x^{{\mathcal{B}}}) \; .
\end{align*}

This shows that~$\tilde{T}_i(\varphi)$ is identically zero
if and only if $H_i(\psi)(x^{{\mathcal{B}}}) = 0$ for
all~$x \in U$. This is equivalent to~$H_i(\psi)$ vanishing
identically on~$V$. Since $T_i \varphi = 0$ if and only if
$\tilde{T}_i \varphi = 0$, we obtain the desired result.
\end{proof}

Parts (\ref{it-1}) and (\ref{it-2}) of Lemma~\ref{lemma:how-to-translate}
have the following consequence.

\begin{corollary}

\label{coro:vector-space-iso}
The map
\[
\begin{array}{rcl}

\{ \mbox{Holomorphic solutions of} \;\;
\Horn({{\mathcal{B}}},c) \; \mbox{on} \; V \}
& \longrightarrow &
\{ \mbox{Holomorphic solutions of} \;
H_{{\mathcal{B}}}(c) \; \mbox{on} \; U \} \\
\psi & \longmapsto & x^c\psi(x^{{\mathcal{B}}})
\end{array} \]
is a vector space isomorphism, that takes Puiseux polynomials to Puiseux polynomials.
\end{corollary}

Finally, we can use the solutions of~$H_A(A\cdot c)$ to construct solutions 
of~$H_{{\mathcal{B}}}(c)$ (and thus of $\Horn({{\mathcal{B}}},c)$). We refer to
\cite[ Section 3]{SST} for background on the canonical series solutions of the
$A$-hypergeometric systems introduced by Gel$'$fand, Kapranov and Zelevinsky.
In the case when~$c$ is generic, these canonical series solutions are fully
supported logarithm-free series.

\begin{theorem}

\label{thm:splitting-series}
 Given a generic parameter vector~$c$,
and $\{ \phi^k : k = 1, \dots ,
\vol(A) \}$ a canonical
basis for the space of solutions of the $A$-hypergeometric
system~$H_A(A\cdot c)$,
there exist linearly
independent, fully supported solutions with disjoint supports
\[ \{ \psi^k_l : k = 1,\dots ,\vol(A) , l = 1, \dots ,g \} \]
of $\Horn({\mathcal{B}},c)$ such that
\[ \phi^{k}=x^c \sum_{l=1}^g \psi^{k}_l (x^{\mathcal{B}})\; ,
\; \mbox{for all} \; \;
k=1,\dots ,\vol(A) \; .\]
Moreover,
no (non trivial) linear combination of the functions~$\psi^k_l$ is ever a
Puiseux polynomial.
This natural decomposition holds
as well for canonical series solutions with logarithms.
\end{theorem}

\begin{proof}

By~\cite[Section 2.5]{SST} and~\cite[Proposition 5.2]{mutsumi2},
a canonical series solution~$\phi$ of the $A$-hypergeometric system~$H_A(A\cdot c)$
is of the form
\begin{equation} \label{eq:phi}
\phi= x^{\alpha} \sum \lambda_{u,v} x^{u} \log(x^v),
\end{equation}
with $A \cdot \alpha = A \cdot c$, and~$v, u \in L=\ker_{\Z}(A)$.
We show that~$\phi$ can
be decomposed as a sum of~$g$ solutions $\psi_1,\dots ,\psi_g$ of
$H_{{\mathcal{B}}}(c)$ such that, if  $\psi_j$, $\psi_l$ are nonzero,
then they have disjoint supports.
Observe that, if $u,v \in L$, then:
\begin{equation}
\label{eqn:goodtheta}
\Big((A\cdot \theta)_j - (A\cdot c)_j\Big) \Big(x^{u+\alpha} \log(x^v)\Big)=0, \quad
\mbox{and}
\end{equation}

\begin{equation}
\label{eqn:no-mixing}
\partial_i  \Big(x^{u+\alpha} \log(x^v)\Big) =
 (u+\alpha)_i \ x^{u+\alpha-e_i} \log(x^v) \, + \, v_i x^{u +\alpha-e_i}.
\end{equation}

Consider the lattice $L_{\mathcal{B}} \subseteq \Z^n$
generated by the columns of~${\mathcal{B}}$,
and its saturation~$L=\ker_{\Z}(A)$,
generated by the columns of a Gale dual~$B$ of~$A$
(that is, the columns of~$B$ form a~$\Z$-basis for the integer
kernel of~$A$). Let $\{ u_l : l=1,\dots ,g \}$ be a system
of representatives for~$L/L_{{\mathcal{B}}}$.
Define
\[\psi_l = x^\alpha \sum_{u \equiv u_l \mbox{\tiny mod} \;L} \lambda_{u,v} x^u \log(x^v).\]
Clearly, $\phi=\psi_1+\cdots +\psi_g$, and the summands have pairwise
disjoint support.
By (\ref{eqn:goodtheta}), each~$\psi_l$ is a solution of the system
of homogeneities $\langle A\cdot \theta - A\cdot c \rangle$. Now we need to check that
each~$\psi_l$ is a solution of the binomial operators $T_1,\dots, T_m$
given by the columns of~${\mathcal{B}}$. Consider
$T_j = \partial^{b^{(j)}_+}-\partial^{b^{(j)}_-}$. Certainly $T_j \phi = 0$.
We apply the operator~$T_j$ to $\phi = \psi_1+\cdots+ \psi_g$, and
observe that terms coming from~$T_j$ applied to~$\psi_l$ cannot
cancel with terms coming from $\partial^{b^{(j)}_+}$ nor from
 $\partial^{b^{(j)}_-}$ applied to~$\psi_{l'}$
if $l \neq l'$.
This is because the exponents of the monomials appearing in
$\Big(\partial^{b^{(j)}_+}\Big)(\psi_l),$ for instance,
are ${b^{(j)}_+}$-translates
of the exponents of the monomials from~$\psi_l$ by (\ref{eqn:no-mixing}),
and $ {b^{(j)}_+} - {b^{(j)}_-} \in L_{{\mathcal{B}}}.$
The lack of cancellation now
follows from the fact that the supports of~$\psi_l$ and~$\psi_{l'}$ are
not congruent modulo~$L_{{\mathcal{B}}}$ by construction.

Now, if we have a canonical
basis $\{ \phi^k, k=1, \dots, {\rm vol}(A)\}$ for the space
of solutions of $H_A (A \cdot c)$ for generic $c \in \C$, they are of the
form
\[\phi^k= x^{\alpha_k} \sum_{u \in  L \cap {\mathcal C}_k} \lambda_{u,v} x^{u}, \]
for different exponents~$\alpha_k$ with respect to a generic weight vector,
and~$u$ ranging over all lattice points in a full dimensional pointed cone
${\mathcal C}_k.$ Notice that, since~$c$ is generic, no pair of the exponents
$\alpha_k$ can differ by an integer vector.
Decompose each $\phi^k = \phi^k_1 + \cdots + \phi^k_g$ as above.
Note that all $\phi^k_l$ are non zero; in fact, the convex hull of all the
supports is full dimensional.
Moreover, the collection $\phi^k_l, \, k=1,\dots, {\rm vol}(A), \, l = 1,\dots, g$
is linearly independent since the supports are disjoint.
By Lemma~\ref{lemma:how-to-translate}, each~$\phi^k_l$ is of the form
$x^c \psi^{k}_l (x^{\mathcal{B}})$, where~$\psi^k_l$ is a solution
of $\Horn({\mathcal{B}},c)$. Clearly, no (non trivial) linear combination of the
functions~$\psi^k_l$ is ever a Puiseux polynomial; in particular, they are
linearly independent.
\end{proof}


\section{Puiseux polynomial solutions of the Horn system and solutions to 
hypergeometric recurrences with finite support}

\label{sec:timur}

Throughout this section we assume that~$m=2.$
Denote by~$\rank_p(J)$ the dimension of the space of Puiseux polynomial
solutions of a $D$-ideal~$J$. 

The first step to compute the dimension of the
space of Puiseux polynomial solutions of $\Horn({\mathcal{B}},c)$
is to observe that such a solution gives rise to a solution of a certain 
system of difference equations. 
A monomial multiple of a Laurent series $\sum_{u\in\Z^m} a(u) y^u$, say 
$y^\gamma \sum_{u\in\Z^m} a(u) y^u$,
is a solution of $\Horn({{\mathcal{B}}},c)$ if and only if its coefficients~$a(u)$
satisfy the recursions:
\begin{equation}
a(u+e_i) Q_i (u+\gamma+e_i) = a(u) P_i (u+\gamma), \,\,\, i=1,\ldots,m.
\label{horn-recurrencies}
\end{equation}
By the support of a solution~$a(u)$ to~(\ref{horn-recurrencies}) we mean the 
set $\{ u: a(u)\neq 0 \}.$ The following proposition is a consequence of 
Proposition~5 in~\cite{pst}. 

\begin{proposition}
Puiseux polynomial solutions of $\Horn({{\mathcal{B}}},c)$ are in one-to-one
correspondence with solutions to~(\ref{horn-recurrencies}) with finite support.
\label{solutions-supports}
\end{proposition}

Let~${\mathcal{B}}[i,j]$ be the
square submatrix of~${\mathcal{B}}$ whose rows are~$b_i$ and~$b_j$, and
let~$c[i,j]$ be the vector in~$\C^2$ whose coordinates are~$c_i$ and~$c_j$. 
We now reduce the computation of the dimension of the
space of Puiseux polynomial solutions to $\Horn({\mathcal{B}},c)$
to the case when~${\mathcal{B}}$ is a $2 \times 2$~matrix.

\begin{lemma}

\label{lemma:ranks-split}
For a generic parameter vector~$c$,
\[
\rank_p(\Horn({\mathcal{B}},c)) = \sum_{i<j} \rank_p(\Horn({\mathcal{B}}[i,j],
c[i,j])).
\]
\end{lemma}

\begin{proof}

We call the support~$S$ of a solution of $\Horn({\mathcal{B}},c)$
{\em irreducible} if there exists no
other solution whose support is a proper nonempty subset of~$S$.
Let~$f(y)$ be a series solution to $\Horn({\mathcal{B}},c)$ with irreducible
support~$S$ and let $s_0 \in S.$
It follows by Theorem~1.3 in~\cite{timur} that if the monomial~$y^{s_0}$
is not present in the series~$f(y)$ then for no $s \in S$ can~$y^s$ be
present in~$f(y).$
This implies that irreducible supports are disjoint.
Indeed, if~$S_1$ and~$S_2$ are irreducible and $s_0\in S_1 \cap S_2$
then there exist solutions~$f_1$ (respectively~$f_2$) of $\Horn({\mathcal{B}},c)$
supported in~$S_1$ (respectively~$S_2$) such that $f_1-f_2$ does not contain~$y^{s_0}.$
But then, since $y^{s_0}$ does not appear in $f_1-f_2$, no monomial in~$S_2$ 
can appear in $f_1-f_2$, and hence $S_1 \backslash S_2$ supports a
solution of $\Horn({{\mathcal{B}}},c)$. This contradicts the fact that~$S_1$ 
was irreducible.

Any Puiseux polynomial solution of $\Horn({\mathcal{B}},c)$ can be written
as a linear combination of polynomial solutions with irreducible supports.
Since Puiseux polynomials with disjoint supports are linearly independent,
it is sufficient to count irreducible supports in order to determine
$\rank_p(\Horn({\mathcal{B}},c))$.

Remember that the equations of the Horn system translate
into recurrence relations~(\ref{horn-recurrencies}) for the coefficients 
of any of its power series solutions.
We refer to~\cite{timur} for a detailed study of these recurrences.
They imply that any coefficient in a  solution of a Horn system
is given by a nonzero multiple of any of its adjacent coefficients, as
long as none of the polynomials~$\bm{P}_i$,~$\bm{Q}_i$ vanish at the
corresponding exponent. This yields that the support of
a  solution must be ``bounded'' by the zeros of these polynomials
in the following sense.
The exponent of a monomial in a  solution must lie
in the zero locus of at least one of the polynomials~$\bm{P}_i$,~$\bm{Q}_i$,
provided that some of the adjacent exponents are not present in the
polynomial solution (See Theorem~1.3 in~\cite{timur}).

Let~$S$ be the support of a Puiseux solution of $\Horn({\mathcal{B}},c)$.
If~$S$ is irreducible, then for a generic vector~$c$ the set~$S$ cannot
meet more than two lines of the form $b_j \cdot \theta_y +c_j -l=0$
corresponding to different parameters~$c_j$.
If it only meets one such line then by Theorem~1.3 in~\cite{timur}
the set~$S$ cannot be finite (in fact, its convex hull is a half-plane
in this case). If~$S$ meets two lines of the above form then all the
other lines can be removed from the picture without affecting
the supports (but not the coefficients) of the Puiseux polynomial
solutions which are generated by this specific pair of lines.
This implies the desired result.
\end{proof}

Now our goal is to compute $\rank_p(\Horn({\mathcal{B}}[i,j],c[i,j]))$.
The first step is to eliminate the cases when this rank is zero.

\begin{lemma}

\label{lemma:supports}
The system $\Horn({\mathcal{B}}[i,j],c[i,j])$ has  non-zero Puiseux polynomial
solutions only if~$b_i$ and~$b_j$ are linearly independent in opposite
open quadrants of~$\Z^2$, or for some special values of $c_i,c_j$ when~$b_i,b_j$ are
linearly dependent and opposite.  The corresponding Puiseux polynomial
solutions of $H_{\mathcal{B}[i,j]}(c[i,j])$ are Taylor polynomials,
that is, polynomials with natural number exponents.
\end{lemma}

\begin{proof}

Corollary~\ref{coro:vector-space-iso} gives a vector space isomorphism
between the solution spaces of the hypergeometric systems
$\Horn({\mathcal{B}}[i,j],c[i,j])$ and
$H_{\mathcal{B}[i,j]}(c[i,j])$ that takes Puiseux polynomials to
Puiseux polynomials.
Thus it is enough to investigate the Puiseux
polynomial solutions of $H_{\mathcal{B}[i,j]}(c[i,j])$. If~$b_i$ and~$b_j$ 
do not lie in the interior of opposite open quadrants, one of the
operators in $H_{\mathcal{B}[i,j]}(c[i,j])$ is of the form
$\partial^{\alpha}-1$ for some $\alpha \in \N^2$. It is clear that
such an operator cannot have a Puiseux polynomial solution.

Now assume that~$b_i$ and~$b_j$ lie in the interior of opposite
quadrants. Let us prove the statement about Taylor polynomials.
We may without loss of
generality assume that~$b_{i1}>0$. If~$b_{i2}<0$, then the change
of variables $\tilde{y}_1 =y_1$,
$\tilde{y_2}=1/y_2$, transforms
$\Horn({\mathcal{B}}[i,j],c[i,j]))$ into a Horn system given by a
$2 \times 2$ matrix whose first row lies in the first open quadrant
of~$\Z^2$. Thus we may assume that $b_{i1},b_{i2}>0$, and consequently
$b_{j1},b_{j2}< 0$, since~$b_i$ and~$b_j$ lie in
opposite open quadrants.

In this case
\[ H_{{\mathcal{B}}[i,j]}(c[i,j]) = \langle \partial_{i}^{b_{i1}}-
\partial_{j}^{-b_{j1}}, \partial_{i}^{b_{i2}}-
\partial_{j}^{-b_{j2}} \rangle ,\]
and this is an ideal in the Weyl algebra with 
generators~$x_i$, $x_j$, $\partial_i$,$\partial_j$.

Let us show that any Puiseux polynomial solution~$f$ of
$H_{{\mathcal{B}}[i,j]}(c[i,j])$ with irreducible support
is actually a Taylor polynomial.
This will imply the statement of the lemma.
Choose $(u_0,v_0) \in \supp(f)$ such that
$\re u_0 = \min \{ \re u : (u,v)\in \supp(f)\backslash \N^2\}$. Then
$(\partial_{i}^{b_{i1}}-\partial_{j}^{-b_{j1}})f$ contains the monomial
$x_i^{u_0-b_{i1}}x_j^{v_0}$ with a nonzero coefficient
unless~$u_0$ is a natural number strictly less than~$b_{i1}$. 
In this case, $v_0 \not \in \N$.
Now, since all the elements of~$\supp(f)$ differ by integer vectors,
and the real part of~$u_0$ is minimal, we have that
$u \in \N$ for all $(u,v) \in \supp(f)\backslash \N^2$.
Now pick
$(u_1,v_1)$ such that the real part of~$v_1$ is minimal, and
conclude that, either~$v_1$ is a natural number strictly
less than~$b_{j1}$ or
$x_i^{u_1}x_j^{v_1-b_{j1}}$ appears with nonzero coefficient in
$(\partial_{i}^{b_{i1}}-\partial_{j}^{-b_{j1}})f=0$.
But now $v \in \N$ for all $(u,v) \in \supp(f)\backslash \N^2$.
We conclude that $\supp(f) \subset \N^2$.

Finally, let us show that if~$b_i$ and~$b_j$ are linearly dependent,
then $\Horn({{\mathcal{B}}[i,j]},c[i,j])$ only has the identically
zero solution as long as~$c$ is generic.
Using the change of variables $\xi_1 = y_{1}^{1/b_{i1}},$ $\xi_2 = y_{2}^{1/b_{i2}},$
we transform the operator~$b_i \cdot \theta_y$ to the operator
$\theta_{\xi_1} + \theta_{\xi_2}.$ By Lemma~\ref{resintheideal}
(to be proved in Section~\ref{sec:holonomic}) there exists a
nonzero polynomial in~$y_1,y_2$ which lies in the ideal
$\Horn({{\mathcal{B}}[i,j]},c[i,j]).$
Thus the only holomorphic solution to the system is the
zero function.
\end{proof}

\begin{example}

\rm
Let us construct the Puiseux polynomial solutions to the system of
equations $\Horn({\mathcal{B}},0),$ where

\[ {\mathcal{B}} =
\left(
\begin{array}{rr}
4  &   5 \\
-3  &  -5 \\
\end{array}
\right).
\]
The system~$H_{{\mathcal{B}}}(0)$ is defined by the operators
\begin{equation}
\frac{\partial^{4}}{\partial x_{1}^{4}} -
\frac{\partial^{3}}{\partial x_{2}^{3}}, \qquad
\frac{\partial^{5}}{\partial x_{1}^{5}} -
\frac{\partial^{5}}{\partial x_{2}^{5}}.
\label{constantcoeff}
\end{equation}

Notice that we may use the parameter~$0$ without loss of
generality. The solutions of $H_{{\mathcal{B}}}(c)$ are exactly the
same as those of~$H_{{\mathcal{B}}}(0)$, and in the case of
$\Horn({\mathcal{B}},c)$, the only effect is a translation of the
supports of the solutions.

The supports of the polynomial solutions to~(\ref{constantcoeff}) are
displayed in Figure~\ref{figure1}.
Two exponents are connected if the
corresponding monomials are contained in a polynomial solution with
irreducible support. Notice that in order to obtain these supports, we
just connected the (empty) circles inside a certain rectangle to other
integer points using the moves given by the columns of~${\mathcal{B}}$.

\begin{figure}[htbp]
\begin{center}

  \psfig{file=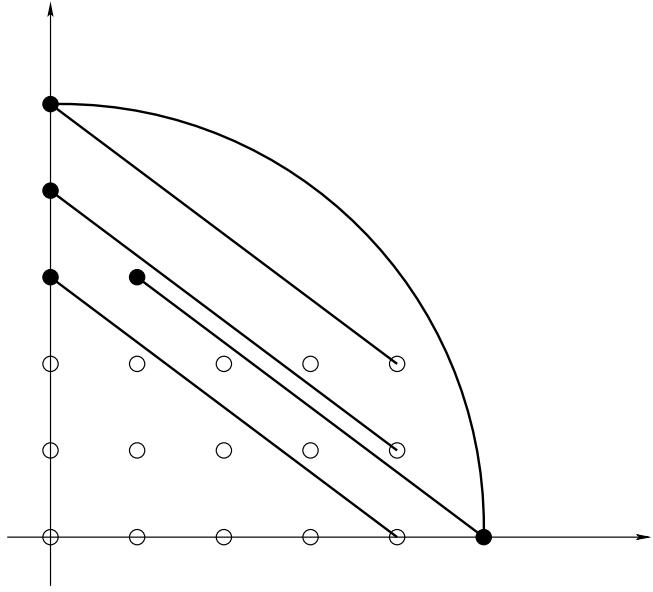}
\end{center}

\caption{The supports of the~$15$ polynomial solutions
to~(\ref{constantcoeff})}

\label{figure1}
\end{figure}

The polynomial solutions to~(\ref{constantcoeff}) are given by
\[1,\quad x_{1}, \quad x_{1}^{2}, \quad x_1^3, \quad x_2, \quad x_1x_2,
\quad x_1^2x_2, \quad x_1^3x_2, \quad x_2^2, \quad x_1x_2^2, \quad
x_1^2 x_2^2, \quad x_1^3x_2^2, \quad x_1^4 + 4x_2^3, \]
\[x_1^4x_2 + x_2^4,\quad  5x_1^4x_2^2 + 2x_1^5 + 2x_2^5 + 40x_1x_2^3.\]

Now let us unravel our isomorphism of solution spaces to obtain the
corresponding solutions of $\Horn({\mathcal{B}},0)$. As in the proof
of the previous lemma, if $\psi = \sum \psi_{\alpha} y^{\alpha}$ is a
Puiseux polynomial solution of  $\Horn({\mathcal{B}},0)$, and
$\psi_{\alpha} \neq 0$, then $\binom{u}{v}={\mathcal{B}}\cdot \alpha \in \N^2$.
But then
\[
\alpha = {\mathcal{B}}^{-1} \cdot \binom{u}{v} = \left(
\begin{array}{rr}
1  &   1 \\
-3/5  &  -4/5 \\
\end{array}
\right) \cdot \binom{u}{v}.
\]
This implies that~$\alpha_1$ is a natural number, and $\alpha_2 \in
(-1/5) \N$. Moreover ${\mathcal{B}} \cdot \alpha \geq 0$. Thus, in
order to find the irreducible supports of the Puiseux polynomial
solutions of $\Horn({\mathcal{B}},0)$, we need to draw the region
 ${\mathcal{B}} \cdot \alpha \geq 0$, plot the points $\alpha \in \N
\times (-1/5) \N$, and connect those points with horizontal and
vertical moves. This is done in Figure~\ref{figure2}. The solid points
belong to the supports of Puiseux polynomials, and the empty circles
and dotted lines correspond to fully supported solutions.
Thus the polynomial solutions to~$\Horn({\mathcal{B}},0)$ are as
follows:

\[1,\quad y_1y_2^{-3/5},\quad y_1^2y_2^{-6/5}, \quad y_1^3y_2^{-9/5},
\quad y_1y_2^{-4/5}, \quad y_1^2y_2^{-7/5}, \quad y_1^3y_2^{-2}, \quad
y_1^4y_2^{-13/5}, \quad y_1^2y_2^{-8/5}, \]
\[ y_1^3y_2^{-11/5}, \quad  y_1^4y_2^{-14/5}, \quad y_1^5y_2^{-17/5}, \quad
y_1^4y_2^{-12/5}+4y_1^3y_2^{-12/5} , \quad
y_1^5y_2^{-16/5}+y_1^4y_2^{-16/5}, \]
\[5 y_1^6y_2^{-4} + 2y_1^5 y_2^{-3} + 2y_1^5 y_2^{-4} + 40y_1^4y_2^{-3}.\]

\begin{figure}[htbp]
\begin{center}
\psfig{file=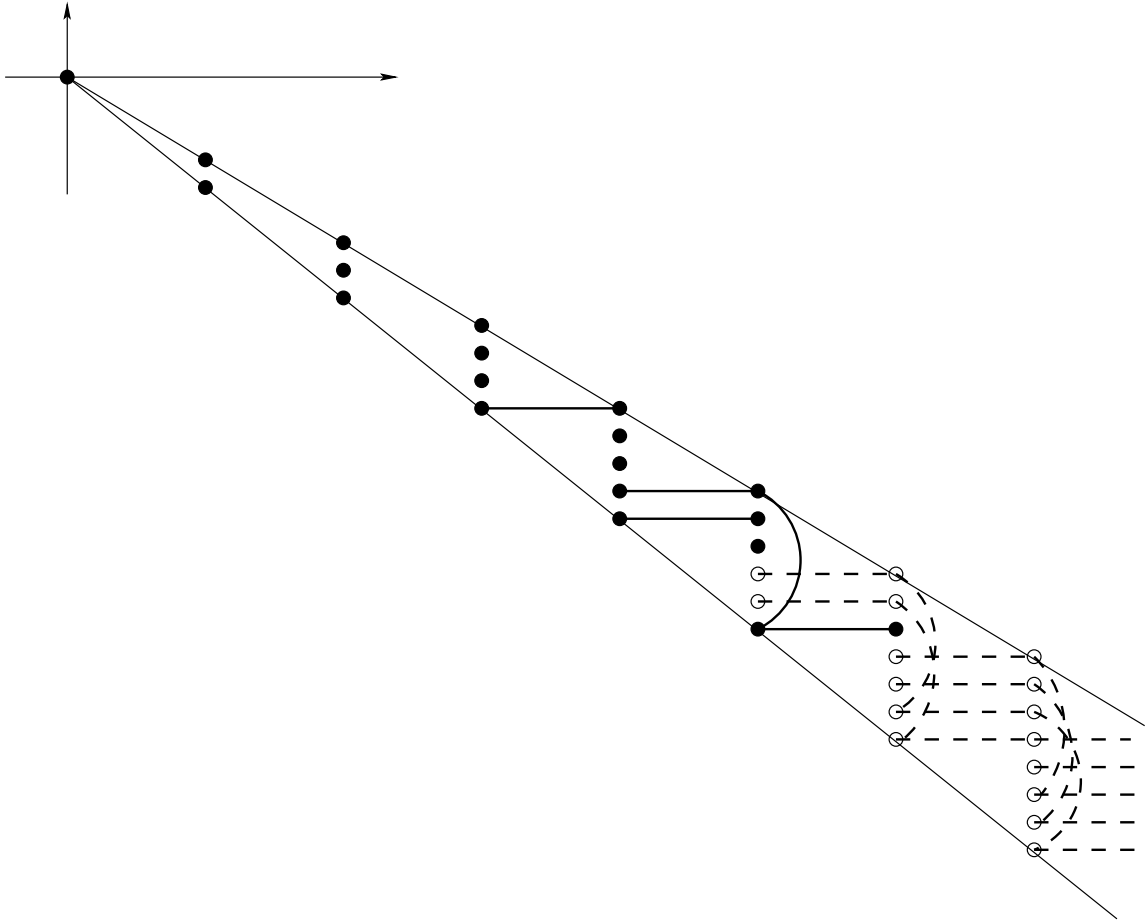}
\end{center}

\caption{The supports of the~15 Puiseux polynomial solutions to $\Horn({\mathcal{B}},0)$
in Example~\ref{ex:puiseux_pols}}

\label{figure2}
\end{figure}

\label{ex:puiseux_pols}
\end{example}

We are now ready to compute
$\rank_{p}(\Horn({{\mathcal{B}}[i,j]},c[i,j]))$.

\begin{lemma}

\label{lemma:rank-for-submatrix}
The dimension of the space of Puiseux polynomial solutions of
the hypergeometric system
$\Horn({{\mathcal{B}}[i,j]},c[i,j])$ equals~$\nu_{ij}$
if the vectors~$b_i$ and~$b_j$ are linearly independent and lie in
opposite open quadrants of~$\Z^2$.
\end{lemma}

\begin{proof}

Suppose that~$b_i$ and~$b_j$ are linearly independent and
lie in opposite open quadrants of~$\Z^2$. As in Lemma~\ref{lemma:supports}, we may assume
that~$b_i$ lies in the interior of the first quadrant (so that~$b_j$
lies in the interior of the third). By Corollary~\ref{coro:vector-space-iso}, it is sufficient
to compute the number of Puiseux polynomial solutions of
$H_{{\mathcal B}[i,j]}(c[i,j]).$

Introduce vectors~$\alpha$,~$\beta$ as follows:
\[ \alpha = \left\{
\begin{array}{ll}
(b_{i1},b_{j1}),  & \;\mbox{if} \; |b_{i1}b_{j2}| > |b_{i2}b_{j1}|, \\
(-b_{i1},-b_{j1}), & \; \mbox{if} \; |b_{i1}b_{j2}| < |b_{i2}b_{j1}|,
\end{array} \right.
\beta = \left\{
\begin{array}{ll}
(-b_{i2},-b_{j2}), & \;\mbox{if} \; |b_{i1}b_{j2}| > |b_{i2}b_{j1}|, \\
(b_{i2},b_{j2}), & \;\mbox{if} \; |b_{i1}b_{j2}| < |b_{i2}b_{j1}|.
\end{array} \right.\]
Furthermore, denote by~${\mathcal{R}}$ the set of points
\[ {\mathcal{R}} = \left\{
\begin{array}{ll}
\{ (u,v) \in \N^2 : u < b_{i2}, \, v< -b_{j1} \} , & \; \mbox{if} \;
 |b_{i1}b_{j2}| > |b_{i2}b_{j1}|, \\
\{ (u,v) \in \N^2 : u < b_{i1}, \, v< -b_{j2} \} , & \; \mbox{if} \;
 |b_{i1}b_{j2}| < |b_{i2}b_{j1}|, \end{array} \right.\]
and call it the {\em base rectangle} of $H_{{\mathcal{B}}[i,j]}(c[i,j])$.
By a {\em path} connecting two points $a, \tilde{a} \in \N^2$ we
mean a sequence $a_1,\dots ,a_k \in \N^2$ such that
$a_1=a$, $a_k=\tilde{a}$ and the difference $a_{i+1}-a_i$ is one of the
vectors $\alpha$, $-\alpha$, $\beta$, $-\beta$. We say that a path
is {\em increasing} if the differences are always one of~$\alpha$,~$\beta$,
and that the path is {\em decreasing} if the differences are always
one of~$-\alpha$,~$-\beta$.
We say that a point in~$\N^2$ is connected with infinity if it can be connected
with another point in~$\N^2$ which is arbitrarily far removed from the origin.

Since the equations defining $H_{{\mathcal{B}}[i,j]}(c[i,j])$ can be
transformed into recurrence relations for the coefficients of a polynomial
solution to this system, it follows that two points can be connected by
a path if and only if the monomials whose exponents are these points
appear simultaneously in a polynomial solution of
$H_{{\mathcal{B}}[i,j]}(c[i,j])$ that has irreducible support.
Notice that if a point in~$\N^2$ is connected with infinity, then
the corresponding monomial cannot be present in any polynomial solution
of~$H_{{\mathcal{B}}[i,j]}(c[i,j])$.

Our next observation is that there are no nonconstant increasing paths
starting at a point of the base rectangle. This can be verified by direct
check of all possible relations between~$|b_{i1}b_{j2}|$,
$|b_{i2}||b_{j1}|$, $b_{i1}$, $b_{i2}$, $b_{j1}$, $b_{j2}$: choosing
the signs of the differences $|b_{i1}b_{j2}| - |b_{i2}||b_{j1}|$,
$b_{i1}-b_{i2}$, $b_{j1}-b_{j2},$ we verify this claim in each of the
eight possible situations.
It follows from this that no two different points in the base rectangle
can be connected by a path, and that no such point is connected with infinity.
Thus, any point in~$\N^2$ is either connected with a unique point in
the base rectangle, or it is connected with infinity. This shows that
the number of polynomial solutions of
$H_{{\mathcal{B}}[i,j]}(c[i,j])$ equals the number of lattice points
in~${\mathcal{R}}$, that is, $\nu_{ij}=\min(|b_{i1}b_{j2}|,|b_{i2}b_{j1}|)$.
\end{proof}

Combining Lemmas~\ref{lemma:ranks-split} and~\ref{lemma:rank-for-submatrix},
we obtain a formula for the dimension of the space of Puiseux
polynomial solutions of $\Horn({\mathcal{B}},c)$.

\begin{theorem}

\label{thm:puiseux}
For a generic parameter~$c$,
\[ \rank_p(\Horn({\mathcal{B}},c))= \sum \nu_{ij} \; ,\]
where the sum runs over pairs of rows~$b_i$ and~$b_j$ of~${\mathcal{B}}$
that are linearly independent and lie in opposite
open quadrants of~$\Z^2$.
\end{theorem}


\section{Solutions of hypergeometric systems arising from lattices}

In this section we consider, for $\beta = A\cdot c$, the lattice hypergeometric system
$I_{\mathcal{B}}+\langle A\cdot \theta - \beta \rangle$. This $D$-ideal is holonomic
for all $\beta \in \C^d$, since its fake characteristic ideal, that is, the ideal
generated by the principal symbols of the generators of~$I_{{\mathcal{B}}}$
and $\langle A\cdot \theta-\beta \rangle$, has dimension~$n$.
In order to compute the holonomic rank of these systems, we need to look at the solutions
of the hypergeometric systems arising from the primary components of~$I_{\mathcal{B}}$.

Let~$\rho$ be a  partial character of~$L/L_{{\mathcal{B}}}$, and let~$I_{\rho}$
be as in Section~\ref{sec:binomialideals}. Define
$H_{\rho}(A\cdot c) =I_{\rho}+\langle A\cdot \theta - A\cdot \beta \rangle$.
In particular, since~$\rho_0$ is the trivial character,
$H_{\rho_0}(A\cdot c)=H_A(A \cdot c)$.

\begin{lemma}

\label{lemma:D-mod-iso}
For $\rho, \rho' \in G_{{\mathcal{B}}}$, the group of partial characters
of $L/L_{{\mathcal{B}}}$, the $D$-modules
$H_{\rho}(\beta)$ and $H_{\rho'}(\beta)$ are isomorphic.
\end{lemma}

\begin{proof}

It is enough to consider the case when $\rho'=\rho_0$, so that $I_{\rho'}=I_{\rho_0} = I_A$.
Given any partial character
$\rho: L \rightarrow \C^*$, let~$p_{\rho}$
be any point in~$X_{\rho}$ all of whose coordinates are nonzero. We define the map
$\tau_{\rho} : D \rightarrow D$ by setting
\[ \tau_{\rho}( \sum x^\alpha \partial^\beta)  \, = \, \sum p_\rho^{\alpha - \beta}  x^\alpha \partial^\beta.\]
It is straightforward to check that~$\tau_{\rho}$ defines an endomorphism of~$D$, 
which is clearly
an isomorphism. It is also easily checked that $\tau_{\rho}(I_A)=I_{\rho}$, and
$\tau_{\rho}(\langle A\cdot \theta - \beta \rangle) = \langle
A\cdot \theta - \beta \rangle$, so that $\tau_{\rho}(H_A(\beta))=H_{\rho}(\beta)$ and 
the $D$-modules~$D/H_A(\beta)$ and~$D/H_{\rho}(\beta)$ are isomorphic.
\end{proof}

\begin{corollary}

\label{coro:reg-holo-Hrho}
If $\rho \in G_{{\mathcal{B}}}$, the $D$-module $D/H_{\rho}(A\cdot c)$
is regular holonomic for all $c \in \C^n$.
\end{corollary}

\begin{proof}

Hotta has shown (see~\cite{equivariant}) that $D/H_A(A\cdot c)$ is regular holonomic
for all parameters $c \in \C^n$, since the condition that the sum of the
rows of~${{\mathcal{B}}}$ equals zero implies that the vector $(1,1,\dots,1)\in \Z^n$
belongs to the row-span of~$A$. Now apply Lemma~\ref{lemma:D-mod-iso}.
\end{proof}

We have shown that the hypergeometric systems arising from the primary
components of the lattice ideal~$I_{\mathcal{B}}$ are regular
holonomic for all parameters. This implies that the solutions of these
systems belong to the Nilsson class~\cite[Ch. 6.4]{jeb1}. We will show that
the solutions of the hypergeometric system
$I_{{\mathcal{B}}}+ \langle A\cdot \theta - \beta \rangle$ satisfy the
same properties.

Recall that  $I_{{\mathcal{B}}}=\cap_{\rho \in G_{\mathcal{B}}} I_{\rho}$,
where $G_{\mathcal{B}}$ is the order~$g$ group of partial
characters, with corresponding ideals~$I_{\rho}$. For any
$\mathcal{J} \subseteq G_{\mathcal{B}}$,
we denote by~$I_{{\mathcal{J}}}$
the intersection $\cap_{\rho \in {\mathcal{J}}} I_{\rho}$.
We first need the following result.

\begin{proposition}

\label{thm:iso}
Let $w \in \N^n \backslash \{ 0 \}$. For  generic~$\beta$, the map
\[ {D}/(I_{\mathcal{J}}+ \langle A\cdot \theta - \beta - A \cdot w \rangle )
\xrightarrow{\quad \cdot \; \partial^w \quad}
{D}/(I_{\mathcal{J}}+ \langle A\cdot \theta - \beta \rangle) \:,
\]
given by right multiplication by~$\partial^w$,
is an isomorphism of left $D$-modules.
\end{proposition}

\begin{proof}

It is sufficient to consider the case when $w = e_i$, so that our map is
right multiplication by~$\partial_i$. In order to use the exact
argument of the proof of~\cite[Theorem 4.5.10]{SST} (the analogous
result for $A$-hypergeometric systems), we need to show that
there exists a nonzero parametric $b$-function (see~\cite[Section
  4.4]{SST}), that is, we need to prove that the
following elimination ideal in the polynomial ring $\C[s_1,\dots,s_d]=\C[s]$
\[ \left( D[s] \; I_{{\mathcal{J}}}+ \langle A \cdot \theta - s
\rangle + D[s] \; \langle \partial_i \rangle \right) \cap \C[s] \]
is nonzero, where~$D[s]$ is the parametric Weyl algebra.
In order to do this, we first go through an intermediate step:
\begin{align*}
\left( D[s] \; I_{{\mathcal{J}}}+ \langle A \cdot \theta - s
\rangle + D[s] \; \langle \partial_i \rangle \right) \cap \C[\theta, s]
& \supseteq \left( D[s](I_{{\mathcal{B}}}+\langle \partial_i\rangle) +
\langle A \cdot \theta - s \rangle \right) \cap \C[\theta, s] \\
& = \left( D[s] (\ini_{-e_i}(I_{\mathcal{B}}+\langle \partial_i
\rangle)) + \langle A \cdot \theta - s \rangle \right) \cap
\C[\theta,s] \\
& = \ini_{(-e_i,e_i,0)} \left(I_{\mathcal{B}}+\langle A\cdot \theta
- s \rangle \right)  \cap \C[\theta, s] \\
& \supseteq \langle [\theta]_{u} : \partial^u \in
\ini_{-e_i}(I_{\mathcal{B}}) \rangle + \langle \theta_i \rangle +
\langle A\cdot \theta - s \rangle \\
& \supseteq \langle [\theta]_{g\  u} : \partial^u \in
\ini_{-e_i}(I_A) \rangle + \langle \theta_i \rangle +
\langle A\cdot \theta - s \rangle \; .\\
\end{align*}
Here $[\theta]_u = \prod_{k=1}^n \prod_{l=0}^{u_k-1} (\theta_k-l)$.
The first containment holds because
$I_{\mathcal{B}} \subseteq I_{\mathcal{J}}$. The next equality is true since
\[ I_{\mathcal{B}} + \langle \partial_i \rangle = \ini_{-e_i}(I_{{\mathcal{B}}}) +
\langle \partial_i \rangle.\]
The equality in the third line holds by the proof of~\cite[Theorem 3.1.3]{SST},
which applies here since~$I_{\mathcal{B}}$ is homogeneous with respect
to the multi-grading given by the columns of~$A$. The next inclusion
is easy to check, given that, for a monomial~$\partial^u$,
$x^u\partial^u = [\theta]_u$. The last containment follows from the fact
that $g \ u \in L_B$ for all $u \in {\rm ker}_\Z (A)$.
Now if we prove that
\[ \left( \langle [\theta]_{g \ u} : \partial^u \in
\ini_{-e_i}(I_A) \rangle + \langle \theta_i \rangle +
\langle A\cdot \theta - s \rangle \right) \cap \C[s] \]
is nonzero, we will be done. But this is a commutative elimination, so all we need to
do is show that the projection of the zero set of $\langle
[\theta]_{g \ u} : \partial^u \in \ini_{-e_i}(I_A) \rangle + \langle \theta_i \rangle +
\langle A\cdot \theta - s \rangle$ onto the $s$-variables is not
surjective.

Observe that the projection of  $\langle [\theta]_{u} : \partial^u \in
\ini_{-e_i}(I_A) \rangle + \langle \theta_i \rangle +
\langle A\cdot \theta - s \rangle$ onto the $s$-variables is not
surjective (by~\cite[Corollary 4.5.9]{SST}). This projection is
clearly the union of affine spaces of different dimensions. But then
the projection that we want is not surjective, since it is obtained
from this one by adding translates of some of the affine spaces
appearing in it. This concludes the proof.
\end{proof}

\begin{theorem}

\label{thm:sols-for-lattice}
For  generic~$\beta$, any solution~$f$ of
$ I_{{\mathcal{J}}}+\langle A\cdot \theta - \beta \rangle$
can be written as a linear combination
\[ f = \sum_{\rho \in {\mathcal{J}}} f_{\rho} \; ,\]
where~$f_{\rho}$ is a solution of $I_{\rho}+\langle A\cdot \theta - \beta \rangle$.
In particular, the solutions of  $I_{\mathcal{B}}+\langle A\cdot \theta - \beta \rangle$
are linear combinations of the solutions of the systems
$I_{\rho}+\langle A\cdot \theta - \beta \rangle$, for $\rho \in G_{{\mathcal{B}}}$.
\end{theorem}

\begin{proof}

We proceed by induction on the cardinality of~${\mathcal{J}}$, the base case
being trivial. Assume that our conclusion is valid for subsets of~$G_{\mathcal{B}}$
of cardinality $r-1 \geq 1$, pick ${\mathcal{J}} \subseteq G_{\mathcal{B}}$ of
cardinality~$r$ and fix $\rho \in {\mathcal{J}}$.

Let~$P$ be an element of
$I_{{\mathcal{J}} \backslash \{ \rho \} }$ such that $P \not \in I_{\rho}$.
Since all of the ideals $I_{\tau}$, $\tau \in G_{\mathcal{B}}$, are homogeneous with respect to
the multi-grading given by~$A$, we may assume that~$P$ is homogeneous, and write
\[P = \lambda_1 \partial^{u^{(1)}} + \cdots +
\lambda_{p-1} \partial^{u^{(p-1)}} + \partial^{w},\]
where
$\lambda_1,\dots ,\lambda_{p-1} \in \C$ and $A\cdot u^{(1)}=A\cdot u^{(2)}\cdots =
A\cdot u^{(p-1)}=A\cdot w$. Notice that the polynomial
\[\bar{P} = \lambda_1 \partial^{u^{(1)}} + \cdots +
\lambda_{p-1} \partial^{u^{(p-1)}}- [ \lambda_1 \rho(u^{(1)}-w)+ \cdots +
\lambda_{p-1} \rho(u^{(p-1)}-w) ]\partial^ w \]
is an element of the ideal~$I_{\rho}$, since
this ideal is generated by all binomials of the form
$\partial^{\alpha} - \rho(\alpha-\gamma) \partial^{\gamma}$, where $A\cdot \alpha = A \cdot \gamma$.
To simplify the notation, set
$-\lambda$ to be the coefficient of $\partial^w$ in $\bar{P}$, that is,
\[ \lambda =  \lambda_1 \rho(u^{(1)}-w)+ \cdots +
\lambda_{p-1} \rho(u^{(p-1)}-w).\]

Now let~$f$ be a solution of $I_{\mathcal{J}}+\langle A\cdot \theta - \beta \rangle$,
and consider the function~$\bar{P} f$. For any
$Q \in I_{{\mathcal{J}}\backslash \{ \rho \}}$, we have
$Q \bar{P} \in I_{{\mathcal{J}}}$. This implies that $Q \bar{P} f = 0$. Furthermore,
noting that~$\bar{P}$ is $A$-homogeneous of multi-degree $A \cdot w$, we conclude
that~$\bar{P} f$ is a solution of $I_{{\mathcal{J}}\backslash \{ \rho \}} +
\langle A \cdot \theta - \beta - A \cdot w \rangle$. Since~$\beta$ is generic, so is
$\beta + A \cdot w$, and by the inductive hypothesis we can write
$\bar{P} f = \sum_{\tau \in {\mathcal{J}} \backslash \{ \rho \} } g_{\tau}$, where 
each~$g_{\tau}$ is a solution of
$I_{\tau} + \langle A \cdot \theta - \beta - A \cdot w \rangle$.

By Proposition~\ref{thm:iso}, $\partial^w$ induces an isomorphism between the solution
spaces of $I_{\tau}+\langle A \cdot \theta - \beta \rangle$ and
$I_{\tau} + \langle A \cdot \theta - \beta - A \cdot w \rangle$, so that we can
find a solution $\tilde{g}_{\tau}$ of $I_{\tau}+\langle A \cdot \theta - \beta \rangle$
such that $\partial^w\tilde{g}_{\tau}=g_{\tau}$. Now
\begin{align*}
\bar{P} \tilde{g}_{\tau} &
= \sum_{i=1}^{p-1} \lambda_i \partial^{u{(i)}} \tilde{g}_{\tau} -
\lambda \partial^w \tilde{g}_{\tau} \\
& = \left(\sum_{i=1}^{p-1} \lambda_i \tau(u^{(i)}-w) - \lambda \right) g_{\tau}.
\end{align*}
The last equality holds because~$\tilde{g}_{\tau}$ is a solution of~$I_{\tau}$, 
and therefore $\partial^{u^{(i)}}-\tau(u^{(i)}-w)\partial^w$ annihilates
it, yielding 
$\partial^{u^{(i)}}\tilde{g}_{\tau} =
\tau(u{(i)}-w)\partial^w \tilde{g}_{\tau} =  \tau(u^{(i)}-w) g_{\tau}$.

Notice that the coefficient $\sum_{i=1}^{p-1} \lambda_i \tau(u^{(i)}-w) - \lambda $
is nonzero, for otherwise we could rewrite~$\bar{P}$ using the sum instead of~$\lambda$,
and conclude that $\bar{P} \in I_{\tau}$. But we know $P \in I_{\tau}$, so
$\bar{P}-P \in I_{\tau}$, a contradiction since this is a nonzero multiple of
$\partial^w$, and the ideal~$I_{\tau}$ contains no monomials. (The fact that
$\bar{P}-P \neq 0$ follows from $\bar{P} \in I_{\rho}$ and
$P \not \in I_{\rho}$).

Finally define $f_{\tau} = \left(\sum_{i=1}^{p-1} \lambda_i \tau(u^{(i)}-w) - \lambda
\right)^{-1} \tilde{g}_{\tau}$, so that~$f_{\tau}$ is a solution of
$I_{\tau} + \langle A \cdot \theta - \beta \rangle$ and
\[ \bar{P} \sum_{\tau \in {\mathcal{J}}\backslash \{ \rho \}} f_{\tau} =
\sum_{\tau \in {\mathcal{J}}\backslash \{ \rho \}} g_{\tau} = \bar{P} f .\]

If $h = f- \sum_{\tau \in {\mathcal{J}}\backslash \{ \rho \}} f_{\tau}$, then~$h$ 
is a solution of $I_{{\mathcal{J}}}+ \langle A \cdot \theta - \beta \rangle$
that satisfies $\bar{P} h = 0$. Now consider~$Ph$. Since
$P\in I_{{\mathcal{J}} \backslash \{ \rho\} }$,  $Ph$ is a solution of
$I_{\rho} + \langle A \cdot \theta - \beta -A\cdot w \rangle$, and a similar argument
as before yields a solution~$f_{\rho}$ of
$I_{\rho} + \langle A \cdot \theta - \beta \rangle$ such that $Ph=Pf_{\rho}$.
Let $\tilde{h} = h- f_{\rho}$, so that $f = \sum f_{\tau} + f_{\rho} +
\tilde{h}$ and $P \tilde{h} = 0$.
But $\bar{P}\tilde{h} = \bar{P} h - \bar{P} f_{\rho} = 0$ since
$\bar{P} \in I_{\rho}$.

Now $P \tilde{h} = \bar{P} \tilde{h}=0$ implies
$(P-\bar{P}) \tilde{h} = 0$, so that $\partial^w \tilde{h} = 0$, because
$P-\bar{P}$ is a nonzero multiple of~$\partial^w$.
But then~$\tilde{h}$ is a solution of
$I_{{\mathcal{J}}}+ \langle A \cdot \theta - \beta\rangle$ that is mapped under
$\partial^w$ to the zero element in the solution space of
$I_{{\mathcal{J}}} + \langle A \cdot \theta - \beta - A\cdot w \rangle$, which, using
the genericity of~$\beta$ and Proposition~\ref{thm:iso}, implies that $\tilde{h}=0$.
Thus we have obtained an expression for~$f$ as a linear combination of solutions
of the systems $I_{\tau}+\langle A \cdot \theta - \beta \rangle$, $\tau \in
\mathcal{J}$, and the proof of the inductive step is finished.
\end{proof}

Considering $\mathcal{J}  = G_{\mathcal{B}},$ we deduce that all solutions
of $D/(I_{{\mathcal{B}}}+\langle A\cdot \theta - \beta \rangle)$ split as a sum
of solutions for each~$I_\rho,$ yielding a kind of converse to 
Theorem~\ref{thm:splitting-series}. We remark that this result is not true 
without the genericity assumption on~$\beta$, since
for certain parameters (for instance for $\beta = 0$, where the constant function~$1$ is
a solution), the solutions to the different
ideals~$H_\rho(\beta)$ are not linearly independent.

\begin{corollary}
\label{coro:lattice-rank}
Suppose that~${\mathcal{B}}$ has zero column sums, and $\beta \in \C^d$ is generic.
Then
\[\rank(I_{\mathcal{B}} + \langle A \cdot \theta - \beta\rangle) \leq g \cdot\vol(A). \]
\end{corollary}

\begin{proof}
Under these hypotheses, the solutions of~$I_{\mathcal{B}}$ are linear combinations
of solutions of the~$g$ systems $I_{\rho}+ \langle A\cdot \theta - \beta\rangle$,
by the previous theorem.
Each of these systems has rank~$\vol(A)$.
\end{proof}


\section{Holonomicity and solutions of the Horn system $H_{{\mathcal{B}}}(c)$}

\label{sec:regular-holo}

In this section we assume that~$m=2$. Our goal is to investigate both the holonomicity
of~$H_{\mathcal{B}}(c)$ and to find out the form of its solutions.
First let us show that~$H_{{\mathcal{B}}}(c)$ is holonomic for
generic~$c$.

\begin{theorem}

\label{thm:hbc-is-holonomic}
Let~$m=2$ and~$c$ generic parameter vector. Then~$H_{{\mathcal{B}}}(c)$
is holonomic.
\end{theorem}

\begin{proof}

Write $I = \langle \partial^{u_+}-\partial^{u_-},
\partial^{v_+}-\partial^{v_-} \rangle$, where~$u$ and~$v$ are the
columns of~${\mathcal{B}}$.
Consider first the case when~${\mathcal{B}}$ has no linearly dependent rows
in opposite open quadrants of~$\Z^2$. Then the ring
\[\frac{\C[x_1,\dots,x_n,z_1,\dots,z_n]}{\langle z^{u_+}-z^{u_-},
z^{v_+}-z^{v_-}\rangle + \langle \sum_{j=1}^n a_{ij}x_jz_j :
i=1,\dots, n-m \rangle}\]
has dimension~$n$ (see Lemma~\ref{lemma:hsop}).
Since the polynomial ring modulo the characteristic ideal
of~$H_{{\mathcal{B}}}(c)$ is a subring of this one, we conclude that
$H_{{\mathcal{B}}}(c)$ is holonomic {\em for all} $c \in \C^m$.

Now assume that~${\mathcal{B}}$ has linearly dependent rows $b_i, b_j$
in opposite open quadrants of~$\Z^2$. In this case, the ideal
$\langle z^{u_+}-z^{u_-},
z^{v_+}-z^{v_-}\rangle + \langle \sum a_{ij}x_jz_j : j=1,\dots ,n-m \rangle$
will have a lower-dimensional component corresponding to the vanishing of~$z_i$ and~$z_j$, 
by the results in Section~\ref{sec:binomialideals} about primary decomposition
of codimension~$2$ lattice basis ideals.

To ensure holonomicity of~$H_{{\mathcal{B}}}(c)$,
we will construct, for each pair~$b_i$,~$b_j$ of linearly dependent rows
of~${\mathcal{B}}$ in opposite open quadrants of~$\Z^2$,
an element of the ideal~$H_{{\mathcal{B}}}(c)$
that contains no~$x_i$,~$x_j$, $\partial_i$, $\partial_j$, and that, for generic~$c$, 
is nonzero. The principal symbol of this element will therefore not depend on~$z_i$ or~$z_j$.

To simplify the notation, assume~$b_1$ and~$b_2$ are linearly dependent in opposite
open quadrants of~$\Z^2$. Then the complementary square submatrix of~$A$ has
determinant zero, so that, by performing row and column operations, we can find
$p,q \in \Q$, $r \in \C$, such that $p \,\theta_1+q\, \theta_2-r$ lies 
in~$H_{{\mathcal{B}}}(c)$. The numbers~$p$ and~$q$ are rational combinations of some
of the elements~$a_{ij}$ of the matrix~$A$, the number~$r$ is a linear combination
of the coordinates of the vector~$c$.

Also, since~$b_1$ and~$b_2$ are linearly dependent, we can find a nonzero element
$w \in L_{\mathcal{B}}$ such that $w_1=w_2=0$. Then we can find two monomials~$m_1$,~$m_2$ 
in $\C[\partial]$ with disjoint supports, that are not divisible by
either~$\partial_1$ or~$\partial_2$ such that
$\partial_1^km_1 (\partial^{w_+}-\partial^{w_-}) \in I$ for some~$k>0$ and
$\partial_2^lm_2 (\partial^{w_+}-\partial^{w_-}) \in I$ for some
$l>0$. This follows from the arguments that proved Proposition~\ref{propo:alicia's-lemma}.
Call $\mu =m_1 (\partial^{w_+}-\partial^{w_-})$ and
$\lambda = m_2 (\partial^{w_+}-\partial^{w_-})$.
Notice that $\mu$, $\lambda$
do not depend on~$\partial_1$,~$\partial_2$.
Then, using $x_1^k\partial_1^k = \theta_1(\theta_1-1)\cdots (\theta_1-k+1)=[\theta_1]_k$
we see that $[\theta_1]_k \mu \in H_{{\mathcal{B}}}(c)$.
Similarly, $[\theta_2]_l \lambda \in  H_{{\mathcal{B}}}(c)$.

Consider the left ideal in the Weyl algebra generated by:
\[ p\, \theta_1+q \,\theta_2-r, [\theta_1]_k \mu, [\theta_2]_l \lambda .\]
This ideal is contained in~$H_{{\mathcal{B}}}(c)$.
Now notice that~$\theta_1$,~$\theta_2$,~$\lambda$ and  $\mu$ are pairwise
commuting elements of~$D_n$. This means that we can think of
$ \langle p\, \theta_1+q\, \theta_2-r, [\theta_1]_k \mu, [\theta_2]_l \lambda
\rangle$ as an ideal in $\C[\theta_1,\theta_2,\partial_3,\dots ,\partial_n]$, which is a
commutative subring of~$D_n$.
We will go one step further and think of~$r$ also as an indeterminate, which
commutes with $\theta_1$, $\theta_2$, $\partial_3,\dots ,\partial_n$.

Finding the element of~$H_{{\mathcal{B}}}(c)$ that we want has now
been reduced to eliminating~$\theta_1$ and~$\theta_2$ from
\begin{equation}
\label{eqn:elimination}
\langle p \, \theta_1+q\, \theta_2-r, [\theta_1]_k \mu, [\theta_2]_l \lambda
\rangle \subset \C[\theta_1,\theta_2,\partial_3,\dots,\partial_n,r] .
\end{equation}
Since the geometric counterpart of elimination is projection,
in order to check that the elimination ideal
\[ \langle  p\, \theta_1+q \, \theta_2-r, [\theta_1]_k \mu, [\theta_2]_l \lambda
\rangle \cap \C[\lambda,\mu,r] \]
is nonzero, we need to show that there exist complex numbers
$\partial_3,\dots,\partial_n$
and~$r$ such that, for all values of $\theta_1, \theta_2 \in \C$,
the tuple $(\theta_1,\theta_2,\partial_3,\dots ,\partial_n,r)$ is not a solution
of (\ref{eqn:elimination}). If $(\partial_3,\dots,\partial_n)$ is generic,
the polynomials~$\mu$ and~$\lambda$ evaluated at that point will be nonzero.
Thus, in order for
$[\theta_1]_k\mu$ to vanish,~$\theta_1$ must be an integer between
$0$ and~$k$. Analogously,~$\theta_2$ must be an integer between~$0$ and~$l$.
But then, for most values of~$r$, $p\, \theta_1+q\, \theta_2-r$ is nonzero.
Thus, the projection of the zero set of (\ref{eqn:elimination}) onto the
$\partial_3,\dots, \partial_n, r$ coordinates is not surjective.
This implies that (\ref{eqn:elimination}) contains an element~$P$ that does not depend
on~$\theta_1$ or~$\theta_2$. Notice that~$P$ does depend (polynomially) on~$r$,
which is itself a linear combination of the coordinates of~$c$. Thus,
for generic~$c$,~$P$ will be nonzero.
Now~$P$ is also an element of the ideal~$H_{{\mathcal{B}}}(c)$,
that does not depend on~$x_1$,~$x_2$, $\partial_1$,~$\partial_2$,
and is nonzero for generic~$c$.
\end{proof}

\begin{example}

Consider the matrix
\[ {\mathcal{B}} = \left( \begin{array}{rr} 1 & 2 \\ -2 & -4 \\ 1 & 1
  \\ 0 & 1 \end{array} \right) .\]
To prove that~$H_{\mathcal{B}}(c)$ is holonomic for generic~$c$, we
need to find an element of~$H_{\mathcal{B}}(c)$ whose principal
symbol does not
vanish if we set $z_1=z_2=0$. To find this element, we follow the
procedure outlined in the proof of the previous theorem. The first
thing we need is an element of $L_{\mathcal{B}}$ with its first two
coordinates equal to zero. The vector $(0,0,-1,1)$ works. It is easy
to check that $\partial_1^2\partial_3^2(\partial_3-\partial_4)$ and
$\partial_2^4(\partial_3-\partial_4)$ are both elements of the lattice
basis ideal~$I$. We can also assume that $(2,1,0,0)$ is a row of the
matrix~$A$. Now what remains is to eliminate~$\theta_1$ and~$\theta_2$
from:
\[ \langle \theta_1(\theta_1-1) \partial_3^2(\partial_3-\partial_4),
\theta_2(\theta_2-1)(\theta_2-2)(\theta_2-3)(\partial_3-\partial_4),
2\theta_1+\theta_2-r \rangle \subseteq
\C[\theta_1,\theta_2,\partial_3,\partial_4,r], \]
where $r=2c_1+c_2$. We perform the elimination on a computer algebra
system to obtain the element:
\[
(2c_1+c_2)(2c_1+c_2-1)(2c_1+c_2-2)(2c_1+c_2-3)(2c_1+c_2-4)(2c_1+c_2-5)
\partial_3^2(\partial_3-\partial_4) \in H_{\mathcal{B}}(c) \]
whose principal symbol
\[
(2c_1+c_2)(2c_1+c_2-1)(2c_1+c_2-2)(2c_1+c_2-3)(2c_1+c_2-4)(2c_1+c_2-5)
z_3^2(z_3-z_4) \]
does not vanish along $z_1=z_2=0$ for generic~$c$.
\end{example}

Our goal now is to characterize all the solutions of the Horn system~$H_{\mathcal{B}}(c)$ 
for generic~$c$. The first step is the following result.

\begin{lemma}

\label{lemma:exact-sequence}
Let~$\alpha$ be as in Proposition~\ref{propo:alicia's-lemma}.
For generic~$c$, the sequence
\begin{equation}
\label{eqn:exact-sequence}
0
\rightarrow
{D}/{(I_{\mathcal{B}}+\langle A \cdot \theta -
  A\cdot (c+\alpha) \rangle)}
\xrightarrow{\; \; \cdot \;\partial^{\alpha}\; \;}
{D}/{H_{{\mathcal{B}}}(c)}
\xrightarrow{\; \; \pi \; \;}
{D}/{(I+\langle \partial^{\alpha}\rangle+ \langle A \cdot \theta -
  A\cdot c \rangle)}
\rightarrow
0,
\end{equation}
where~$\pi$ is the natural projection, is exact.
\end{lemma}

\begin{proof}

The only part of exactness that is not clear is that
right multiplication by~$\partial^{\alpha}$ is injective (it is well
defined since $\partial^{\alpha} I_{{\mathcal{B}}} \subseteq I$).
To see this, consider the following commutative diagram:
\diagram[PostScript=dvips]
 0 & \rTo &
{D}/ (I_{{\mathcal{B}}}+\langle  A\cdot \theta - A \cdot (c+\alpha) \rangle)
& \rTo^{\cdot \, \partial^{\alpha}} &
{D} /(I_{{\mathcal{B}}}+\langle  A\cdot \theta - A \cdot c \rangle) &
 \rTo 0 \\
 &  &  & \rdTo_{\cdot \, \partial^{\alpha}} & \dTo  &  \\
& & &  & {D} /{H_{\mathcal{B}}(c)}  &
\enddiagram
where the vertical
arrow is the natural inclusion. The upper row
of the diagram is exact by Theorem~\ref{thm:iso}, since~$c$ is
generic. But then the commutativity implies that the 
diagonal arrow is injective.
\end{proof}

\begin{lemma}

\label{lemma:laura's-lemma}
Let $u,v \in \N^n$ such that $\langle \partial^u,\partial^v \rangle$
is a complete intersection. If~$c$ is generic, then
\[ \langle \partial^u,\partial^v \rangle + \langle A \cdot \theta - A
\cdot c \rangle \]
is a holonomic system of differential equations, whose
solution space has a basis of Puiseux monomials.
\end{lemma}

\begin{proof}
It is enough to show that the system
\[ \langle x^u\partial^u,x^v\partial^v \rangle + \langle A \cdot \theta - A
\cdot c \rangle \]
satisfies the desired properties since~$x^u$ and~$ x^v$ are units in~$\C(x)$.

Now
\[ \langle x^u\partial^u,x^v\partial^v \rangle + \langle A \cdot \theta - A
\cdot c \rangle
=  \langle [\theta]_u, [\theta]_v \rangle + \langle A \cdot \theta - A
\cdot c \rangle = D \cdot F,\]
where
\[ [\theta]_u = \prod_{k=1}^n \prod_{l=0}^{u_k-1} (\theta_k-l) \; , \]
and
\[ F = \langle [\theta]_u, [\theta]_v \rangle + \langle A \cdot \theta - A
\cdot c \rangle \subseteq \C[\theta] .\]
This means that~$D\cdot F$ is a Frobenius ideal (see~\cite[Section
  2.3]{SST}). By~\cite[Proposition 2.3.6, Theorem 2.3.11]{SST}, if we
can show that~$F$ is artinian and radical, it will follow that~$D \cdot F$
is holonomic, with solution space spanned by $\{ x^p : p \in
{\mathcal{V}}(F)\}$, where~${\mathcal{V}}(F)$ is the zero set of the
ideal $F \subseteq \C[\theta]$, and we will be done.

To show that~$F$ is artinian and radical, we proceed as in~\cite[Theorem 3.2.10]{SST}. 
Let $p \in {\mathcal{V}}(F)$. Then there
exist $1 \leq i < j \leq n$ such that~$p_i$ and~$p_j$ are nonnegative integers
between zero and $\max\{u_i,v_i\}$, $\max\{u_j,v_j\}$
respectively. This follows from $[\theta]_{u}(p)=[\theta]_{v}(p)=0$
and  the fact that~$u$ and~$v$ have disjoint supports, because
$\langle\partial^u, \partial^v \rangle$ is a complete
intersection. Since~$c$ is generic, the minor of~$A$ complementary to~$\{ i,j\}$ 
must be nonzero (otherwise the equations $\theta_i=p_i$,
$\theta_j=p_j$ and $A\cdot \theta = A\cdot c$ would be
incompatible). Hence its $i$-th and $j$-th coordinates determine~$p$
uniquely in~${\mathcal{V}}(F)$.
\end{proof}

\begin{remark}
If all maximal minors of~$A$ are nonzero, the above lemma holds
without restriction on~$c$.
\end{remark}

\begin{theorem}
\label{theorem:system+monomial}
Write $I=\langle \partial^{u_+}-\partial^{u_-}
,\partial^{v_+}-\partial^{v_-} \rangle$, where~$u$ and~$v$ are the
columns of~${\mathcal{B}}$. Let~$\partial^{\alpha}$ be a monomial satisfying:
\begin{equation}
\label{eqn:alpha}
\alpha_i > 0 \Longrightarrow u_i > 0.
\end{equation}
Then, for generic~$c$, the $D$-ideal $I+\langle \partial^{\alpha}
\rangle + \langle A\cdot \theta - A \cdot c \rangle$ has only Puiseux
polynomial solutions.
\end{theorem}

\begin{proof}
We proceed by induction on $|\alpha|=\alpha_1+\cdots+\alpha_n$, the
length of~$\alpha$.
If $|\alpha| \leq \min\{u_i:u_i>0\}$, in particular, if $|\alpha|=1$
(recall that $|u|=0$), then $\partial^{\alpha}$ divides
$\partial^{u_+}$, so that all solutions of $I+\langle \partial^{\alpha}
\rangle + \langle A\cdot \theta - A \cdot c \rangle$ are solutions of
$\langle \partial^{\alpha},\partial^{u_-} \rangle + \langle A\cdot
\theta - A \cdot c \rangle$. But the latter ideal has only Puiseux
polynomial solutions by Lemma~\ref{lemma:laura's-lemma}, since~$c$ is generic.

Assume now that our result is true for length~$s$ and let~$\alpha$ be
of length $s+1$ satisfying (\ref{eqn:alpha}). Choose~$i$ such that
$\alpha_i>0$ (and so $u_i>0$), and let~$\varphi$ be a solution of $I+\langle \partial^{\alpha}
\rangle + \langle A\cdot \theta - A \cdot c \rangle$. The function
$\partial_i \varphi$ is a solution of $I+\langle \partial^{\alpha-e_i}
\rangle + \langle A\cdot \theta - A \cdot c -A\cdot e_i \rangle$. But
$|\alpha-e_i|=s$ and $c+e_i$ is still generic, so the inductive
hypothesis implies that $\partial_i \varphi$ is a Puiseux
polynomial. Write:
\[ \partial_i \varphi = \sum_{l=0}^{N_0} g_l^{(0)} x_i^l +\sum_{l=0}^{N_1}
g_l^{(1)} x_i^{\mu_1+l} + \cdots +\sum_{l=0}^{N_t} g_l^{(t)} x_i^{\mu_t+l} ,\]
where the $g_l^{(k)}$ are Puiseux polynomials, constant with respect to~$x_i$, $t$ 
is a natural number, and $\mu_1,\dots,\mu_t \in \C$ are
nonintegers with noninteger pairwise differences.
Then
\begin{equation}
\label{eqn:expr-for-varphi}
\varphi = \sum_{l=0}^{N_0} g_l^{(0)} \frac{x_i^{l+1}}{l+1} +\sum_{l=0}^{N_1}
g_l^{(1)} \frac{x_i^{\mu_1+l+1}}{\mu_1+l+1} + \cdots +\sum_{l=0}^{N_t} g_l^{(t)}
\frac{x_i^{\mu_t+l+1}}{\mu_t+l+1} + G(x_1,\dots,\hat{x_i},\dots,x_n).
\end{equation}

If we prove that~$G$ is a Puiseux polynomial, it will follow that so
is~$\varphi$, and the proof will be finished.
We know that~$\varphi$ is a solution of $\langle A\cdot \theta - A
\cdot c \rangle$. By construction, so is $\varphi - G$. Then~$G$ is a
solution of  $\langle A\cdot \theta - A \cdot c \rangle$. Recall
that $\partial_i G=0$.

We also know that $\partial^{u_+} \varphi = \partial^{u_-} \varphi$.
We want to compare the coefficients of the integer powers of~$x_i$ in
the expressions we obtain by applying~$\partial^{u_+}$
and~$\partial^{u_-}$ to (\ref{eqn:expr-for-varphi}). Since we are only
looking at the integer powers of~$x_i$, we need only look at
$\sum_{l=0}^{N_0} g^{(0)}_l(x_i^{l+1}/(l+1)) +G$.
\begin{equation}
\label{eqn:expr1}
\partial^{u_+} 
\left(
\sum_{l=0}^{N_0} g^{(0)}_l\frac{x_i^{l+1}}{l+1} +G
\right) = 
\sum_{l=0}^{N_0} l(l-1)\cdots (l+2-u_i)
(\partial^{u_+ -u_ie_i} g^{(0)}_l) x_i^{l+1-u_i} .
\end{equation}

Notice that there is no~$G$ in the above expression, since
$\partial_i G=0$ and $u_i >0$. Also, the highest power of~$x_i$
appearing in (\ref{eqn:expr1}) is $x_i^{N_0+1-u_i}$.

\begin{equation}
\label{eqn:expr2}
\partial^{u_-} \left(\sum_{l=0}^{N_0} g^{(0)}_l\frac{x_i^{l+1}}{l+1} +G
\right) = \sum_{l=0}^{N_0} (\partial^{u_-}g^{(0)}_l)
\frac{x_i^{l+1}}{l+1} + \partial^{u_-}G .
\end{equation}

We equate the coefficients of $x_i^{l+1}$ in (\ref{eqn:expr1}) and
(\ref{eqn:expr2}) to obtain:
\begin{equation}
\label{eqn:expr3}
\frac{\partial^{u_-}g^{(0)}_l}{l+1} = (l+u_i)\cdots (l+2)
\partial^{u_+-u_ie_i} g^{(0)}_{l+u_i} \; , \; \mbox{for} \; l=0,\dots ,N_0-u_i.
\end{equation}
If $l=N_0+1-u_i,\dots,N_0$, then $\partial^{u_-}g^{(0)}_l=0$. Also,
\[\partial^{u_-} G = (u_i-1)(u_i-2)\cdots 2 \cdot 1 \cdot
\partial^{u_+-u_ie_i}g^{(0)}_{u_i-1}.\]
Applying $\partial^{u_-}$ to (\ref{eqn:expr3}), we see that
\[ \partial^{2u_-} g^{(0)}_l = (l+u_1)\cdots(l+2)(l+1)
\partial^{u_+-u_ie_i}\partial^{u_-}g^{(0)}_{l+u_i} = 0 \; , \;
\mbox{for}\; l = N_0+1-2u_i,\dots,N_0-u_i.\]
Applying $\partial^{u_-}$ enough times, we conclude that,
if $ku_i>N_0+1,$ then $\partial^{ku_-}G=0.$
But now,~$G$ is a solution of $\langle \partial_i,\partial^{ku_-}
\rangle + \langle A \cdot \theta - A \cdot c \rangle$, and~$c$ is generic. By Lemma
\ref{lemma:laura's-lemma},~$G$ is a Puiseux polynomial.
\end{proof}

\begin{proposition}
\label{propo:solution-decomposition}
Let $\alpha$ be as in Proposition \ref{propo:alicia's-lemma} (in particular,
$\alpha$ satisfies (\ref{eqn:alpha})), let $c$ be generic, and let $f$
be a solution of $H_{\mathcal{B}}(c)$. Then $f=g+h$, where $g$ is a solution of the
lattice hypergeometric system $I_{\mathcal{B}}+\langle A\cdot \theta - A\cdot c \rangle$
and $h$ is a solution of $I+\langle \partial^{\alpha} \rangle +
\langle A\cdot \theta - A\cdot c \rangle$.
\end{proposition}

\begin{proof}
Let $\psi = \partial^{\alpha} f$. Then $\psi$ is a solution of
$I_{B}+\langle A\cdot \theta - A\cdot (c+\alpha) \rangle$. This is because
the $D$-module map
\[ \frac{D}{I_{\mathcal{B}}+\langle A\cdot (c+\alpha) \rangle}
\xrightarrow{\quad \cdot \;\partial^{\alpha} \quad } \frac{D}{H_{\mathcal{B}}(c)} \]
induces a vector space map between the solution spaces of $H_{\mathcal{B}}(c)$ and
$I_{\mathcal{B}}+\langle A\cdot \theta - A\cdot (c+\alpha) \rangle$.

Now by Lemma \ref{lemma:D-mod-iso}, right multiplication by $\partial^{\alpha}$
is an $D$-module isomorphism between
$D/(I_{\mathcal{B}}+ \langle A \cdot \theta - A\cdot (c+\alpha) \rangle)$
and
$D/(I_{\mathcal{B}}+ \langle A \cdot \theta - A\cdot c \rangle)$, 
so there exists $Q \in D$ and
$P \in I_{\mathcal{B}} + \langle A \cdot \theta - A\cdot (c+\alpha) \rangle$ such that 
$\partial^{\alpha} Q = 1 + P$. Let $g = Q \psi$.
Then $g$ is a solution of $I_{\mathcal{B}} + \langle A\cdot \theta - A\cdot c \rangle$,
and
\begin{equation}
\label{eqn:dalpha}
\partial^{\alpha} g = \partial^{\alpha} Q \psi = (1+P) \psi = \psi
=\partial^{\alpha} f
\end{equation}
where the next to last equality holds because $P \in I_{\mathcal{B}} +
\langle A\cdot \theta -A\cdot (c+\alpha) \rangle$. Now let $h = f-g$. All we need to
finish this proof is to show that $h$ is a solution of
$I + \langle \partial^{\alpha} \rangle + \langle A\cdot \theta - A \cdot c \rangle$.
But, since $I \subset I_{\mathcal B}$, $g$ is also a solution of $H_{\mathcal{B}}(c)$, 
and thus so is $h$. Moreover $\partial^{\alpha} h =0$ by (\ref{eqn:dalpha}).
\end{proof}

\begin{corollary}
\label{coro:horn-rank-bound}
For generic~$c$, we have
\[ \rank(H_{\mathcal{B}}(c)) \leq g \cdot \vol(A) + \sum \nu_{ij}, \]
where the sum runs over pairs of linearly independent rows of~${\mathcal{B}}$
in opposite open quadrants of~$\Z^2$.
\end{corollary}

\begin{proof}
By Proposition~\ref{propo:solution-decomposition}, the solution space of
$H_{\mathcal{B}}(c)$ is contained in the sum of the solution spaces of
$I_{\mathcal{B}}+ \langle A\cdot \theta - A\cdot c\rangle$ and
$I+ \langle \partial^{\alpha} \rangle+ \langle A\cdot \theta - A\cdot c\rangle$.
The first solution space has rank at most
$g \cdot \vol(A)$ by Corollary~\ref{coro:lattice-rank}.
The second solution space contains only Puiseux polynomials and
therefore has rank at most $\rank_p(H_{\mathcal{B}})=\sum \nu_{ij}$ by Theorem
\ref{thm:puiseux}.
\end{proof}


\section{Initial ideals, indicial ideals and holonomic ranks}

\label{sec:bound}

In this section we finish the proofs of our rank formulas for generic parameters,
by showing the reverse inequalities in Corollaries~\ref{coro:lattice-rank} and
\ref{coro:horn-rank-bound}.
We will assume $m=2$ when dealing with Horn systems, although the arguments will work
for general~$m$ as long as~$I$ is a complete intersection and~$H_{{\mathcal{B}}}(c)$
is holonomic for generic~$c$.

Our main tool will be the fact that holonomic rank is lower semicontinuous when
we pass to initial ideals with respect to weight vectors of the form $(-w,w)$;
this is~\cite[Theorem 2.2.1]{SST}.
For an introduction to initial ideals in the Weyl algebra, including algorithms, see
\cite[Chapters 1 and 2]{SST}.

\begin{theorem}[Theorem 2.2.1~\cite{SST}]
\label{thm:semic-holo-rank}
If~$J$ is a holonomic $D_n$-ideal, and~$w$ is a generic weight vector, then the
initial $D_n$-ideal $\ini_{(-w,w)}(J)$ is also holonomic, and
\[ \rank(\ini_{(-w,w)} (J )) \leq \rank(J) .\]
\end{theorem}

\begin{remark}
If we assume that~$J$ is {\em regular holonomic}, then equality will hold in the above
theorem.
\end{remark}

Our goal is now to compute the holonomic ranks of $\ini_{(-w,w)}(H_{{\mathcal{B}}}(c))$
and $\ini_{(-w,w)}(I_{{\mathcal{B}}})+\langle A\cdot \theta - A\cdot c \rangle$ for
generic~$c$. In order to do this, we introduce indicial ideals, which are modifications
of initial ideals, and have the advantage of belonging to
the (commutative) polynomial ring~$\C[\theta]$.

\begin{definition}
If~$J$ is a holonomic left $D_n$-ideal, and~$w$ is a generic
weight vector, the {\em indicial ideal} of~$J$ is
\[ \ind_w(J)= R \cdot \ini_{(-w,w)}(J) \cap \C[\theta_1,\dots,\theta_n], \]
where~$R$ is the ring of linear partial differential equations with
rational function coefficients.
\end{definition}

A $D_n$-ideal whose generators belong to $\C[\theta]=\C[\theta_1,\dots ,\theta_n]$
is called a {\em Frobenius ideal}. The commutative ideal in~$\C[\theta]$ given by the
generators of a Frobenius ideal is called the {\em underlying commutative ideal}.
The following theorem justifies our interest in indicial ideals.

\begin{theorem}[Theorem 2.3.9~\cite{SST}]
Let~$J$ be a holonomic $D_n$-ideal and~$w$ a generic weight vector. Then
$D_n \cdot \ind_w(J)$ is a holonomic Frobenius ideal whose rank equals
$\rank(\ini_{(-w,w)}(J))$.
\end{theorem}

Finally, computing the rank of a holonomic Frobenius ideal (such as $\ind_w(J)$
for holonomic~$J$) is a commutative operation.

\begin{proposition}[Proposition 2.3.6~\cite{SST}]
Let $D_nF$ be a Frobenius ideal, where $F \subset \C[\theta]$
is the underlying commutative ideal. Then~$D_nF$ is holonomic
if and only if~$F$ is zero-dimensional, in which case
\[ \rank(D_nF) = \deg(F). \]
\end{proposition}

Although indicial ideals are extremely useful, they
are hard to get a hold of in general.
However, for generic parameters, we know explicitly what the indicial
ideal of an $A$-hypergeometric system is (\cite[Corollary  3.1.6]{SST}),
and the same ideas work for the case of Horn systems
and hypergeometric systems arising from lattices.

\begin{theorem}
For generic parameters~$c$, we have
\[ \ind_w(H_{{\mathcal{B}}}(c)) = \big( (R \cdot \ini_w(I)) \cap \C[\theta] \big)+
\langle A \cdot \theta -
A \cdot c \rangle, \]
and
\[ \ind_w(I_{\mathcal{B}}+\langle A\cdot \theta - A \cdot c \rangle) =
\big( (R \cdot \ini_w(I_{\mathcal{B}})) \cap \C[\theta] \big)+
\langle A \cdot \theta -
A \cdot c \rangle. \]
\end{theorem}

\begin{proof}
The proof of the analogous fact for $A$-hypergeometric systems follows
from~\cite[Theorem 3.1.3 and Proposition 3.1.5]{SST}. But~\cite[Proposition 3.1.5]{SST} 
carries over to the cases that interest
us without any modification in its proof. Moreover
the proof of~\cite[Theorem 3.1.3]{SST} only uses the fact that~$I_A$
is homogeneous with respect to the multi-grading given by the columns
of~$A$, a property that both~$I$ and~$I_{\mathcal{B}}$ satisfy.
\end{proof}

Our next goal is to compute the primary decomposition of
the indicial ideals of~$H_{\mathcal{B}}(c)$ and $I_{\mathcal{B}}+\langle A\cdot
\theta - A\cdot c \rangle$
when~$c$ is generic.
The first step is to recall the definition of
certain combinatorial objects that
correspond to the irreducible components of a monomial ideal in a polynomial ring.

\begin{definition}
Let~$M$ be a monomial ideal in the polynomial ring $\C[\partial_1, \dots
,\partial_n]$. A {\em standard pair} of~$M$ is a pair
$(\partial^{\eta}, \sigma)$, where~$\sigma$ is a possibly empty
subset of $\{ 1, \dots ,n \}$,
that satisfies
\begin{enumerate}
\item $\eta_i = 0$ for all $i \in \sigma$;
\item for any choice of integers~$\mu_j \geq 0$,
$j \in \sigma$, the monomial
$\partial^{\eta} \prod_{j \in \sigma} \partial_j^{\mu_j}$ is not in~$M$;
\item for all $l \not \in \sigma$, there exist integers
$\mu_l \geq 0$ and $\mu_j \geq 0$, $j \in \sigma$, such that
$\partial^{\eta}\partial_l^{\mu_l}\prod_{j \in \sigma} \partial_j^{\mu_j}$
lies in~$M$.
\end{enumerate}
\end{definition}

We denote the set of standard pairs of a monomial ideal~$M$ by~$S(M)$.
By~\cite[Equation (3.2)]{StuTruVo},
\[ M =\bigcap_{(\partial^{\eta},\sigma) \in S(M)} \langle \partial_i^{\eta_i+1}
: i \not \in \sigma \rangle . \]
The prime ideal $\langle \partial_i : i \not \in \sigma \rangle$ is associated
to~$M$ if and only if there exists a standard pair of the form
$(\cdot , \sigma)$ in~$S(M)$.
A standard pair $(\partial^{\eta},\sigma)$ is called {\em top dimensional}
if $\langle \partial_i : i \not \in \sigma \rangle$ is
a minimal associated prime of~$M$, it is called {\em embedded} otherwise.
It is clear from the above formula that the degree of~$M$ is equal to the
cardinality of the set of top dimensional standard pairs of~$M$.

Now, since the ideals~$I$ and~$I_{\mathcal{B}}$ are unmixed ($I$ is a
complete intersection, and the associated primes of~$I_{\mathcal{B}}$
are all isomorphic to~$I_A$),
all of the minimal primes of all the initial ideals
of~$I$ have the same dimension,~$d$ (see~\cite[Corollary 1]{KalStu}),
and the same holds for~$I_{\mathcal{B}}$.
This means that a standard pair~$(\partial^{\eta}, \sigma)$
of either $\ini_w(I)$ or $\ini_w(I_{\mathcal{B}})$ is top dimensional
if and only if~$\# \sigma = d$.

Let~$T(\ini_w(I))$ be the set of top dimensional
standard pairs $(\partial^{\eta},\sigma)$
of~$\ini_w(I)$ such that the rows of~${{\mathcal{B}}}$ indexed by
$i \not \in \sigma$ are linearly independent.

Notice that if $(\partial^{\eta},\sigma)$ is a top-dimensional
standard pair of $\ini_w(I_{\mathcal{B}})$, then the rows of~${\mathcal{B}}$ 
indexed by $i \not \in \sigma$ are linearly
independent (the proof of~\cite[Lemma 2.3]{shared} works for lattice
ideals too). Then $T(\ini_w(I_{\mathcal{B}}))$ equals the set of
top-dimensional standard pairs of $\ini_w(I_{\mathcal{B}})$.

Given a standard pair in either $T(\ini_w(I))$ or $T(\ini_w(I_{\mathcal{B}}))$,
and an arbitrary parameter vector~$c$, there exists a unique vector~$v$ such that 
$A \cdot v = A \cdot c$, and $v_k=\eta_k$, $v_l=\eta_l$.

Suppose that $(\partial^{\eta}, \sigma)$ is a standard pair of
$\ini_w(I)$ that does not belong to~$T(\ini_w(I))$. Then either
$\# \sigma < m$ or $\# \sigma = n-2$ and the columns of
${\mathcal{B}}$
corresponding to the indices not in~$\sigma$ are linearly dependent.
In both of these cases, for a generic choice of~$c$, the system
$A \cdot v = A \cdot c$, $v_i = \eta_i$ for $i \not \in \sigma$,
has no solutions. The same holds for standard pairs not in 
$T(\ini_w(I_{\mathcal{B}}))$.

We can now describe the primary decomposition of the indicial ideals
of~$H_{{\mathcal{B}}}(c)$ and $I_{\mathcal{B}}+\langle A \cdot \theta
- A \cdot c\rangle$
with respect to~$w$, in analogy to~\cite[Theorem 3.2.10]{SST}.

\begin{proposition}
For a generic parameter~$c$,
the indicial ideal of~$H_{{\mathcal{B}}}(c)$
with respect to~$w$ equals the
following intersection of maximal ideals:
\begin{equation}
\label{eqn:fakeind}
\bigcap_{(\partial^{\eta},\sigma) \in T(\ini_w(I))}
\big (\langle \theta_i - \eta_i : i \not \in \sigma \rangle +
\langle A \cdot \theta - A \cdot c \rangle \big) ,
\end{equation}
and the indicial ideal of $I_{\mathcal{B}}+\langle A\cdot \theta -
A\cdot c\rangle$ equals:
\begin{equation}
\label{eqn:fakeind2}
\bigcap_{(\partial^{\eta},\sigma) \in T(\ini_w(I_{\mathcal{B}}))}
\big (\langle \theta_i - \eta_i : i \not \in \sigma \rangle +
\langle A \cdot \theta - A \cdot c \rangle \big) .
\end{equation}

\end{proposition}

\begin{proof}
We prove the statement for the indicial ideal of~$H_{\mathcal{B}}(c)$. 
The other indicial ideal is computed in exactly in the same manner.

By~\cite[Corollary 3.2.3]{SST}, the indicial ideal is
\[ J = \langle A \cdot \theta - A \cdot c \rangle +
\bigcap_{(\partial^{\eta},\sigma) \in S(\ini_w(I))} \langle
\theta_i - \eta_i: i \not \in \sigma \rangle .\]
It is clear that the ideal (\ref{eqn:fakeind}) is radical.
If we show that it has the same zero set as~$J$,
and that~$J$ has no multiple roots, we will be done.

Let~$v$ be a zero of~$J$. Then $A \cdot v = A \cdot c$, and
for some $(\partial^{\eta}, \sigma) \in S(\ini_w(I))$, we have that
$v_i = \eta_i$ for all $i \not \in \sigma$. Since our
parameter~$c$ is generic, we must have that $(\partial^{\eta},\sigma)$
belongs to~$T(\ini_w(I))$. These are exactly the roots of the
ideal (\ref{eqn:fakeind}).
It also follows from the genericity of~$c$ that all the zeros of~$J$ 
are distinct, and the proof is finished.
\end{proof}

Notice that the degree of~$\ini_w(I)$ is~$d_1 \cdot d_2$,
since it coincides with the degree of the complete intersection~$I$.
Then the cardinality of the set of top dimensional standard pairs is
exactly~$d_1 \cdot d_2$.
This and the previous proposition imply the following result.

\begin{corollary}
\label{coro:first-bound}
Let~$\nu$ be the sum of the multiplicities of the minimal primes
of~$\ini_w(I)$ corresponding to linearly dependent sets of
two rows of~${\mathcal{B}}$.
For a generic parameter vector~$c$, the degree of the fake indicial ideal is exactly
$d_1 \cdot d_2 - \nu$. Therefore,
\[ \rank(H_{{\mathcal{B}}}(c))= \rank(\Horn({\mathcal{B}},c)) =
d_1\cdot d_2-\nu = \#T(\ini_w(I)) .\]
\end{corollary}

Our desired formula for the generic rank of a bivariate Horn system now follows from
Proposition~\ref{propo:lattice-basis-description}.

\begin{theorem}

\label{thm:rank-formula}
For generic~$c$ and~$m=2$,
\[ \rank(H_{{\mathcal{B}}}(c))=\rank(\Horn({{\mathcal{B}}},c))=
d_1 \cdot d_2 - \sum \nu_{ij} \;,\]
where the sum runs over linearly dependent rows of~${{\mathcal{B}}}$ that
lie in opposite open quadrants of~$\Z^2$.
\end{theorem}

\begin{proof}

By Proposition~\ref{propo:lattice-basis-description}, the sum of the multiplicities
of the minimal primes of~$I$ corresponding to linearly dependent rows 
of~${{\mathcal{B}}}$ is the sum of the corresponding
indices $\sum \nu_{ij}$. This implies that
\[\deg(\ind_w(H_{{\mathcal{B}}}(c)) = d_1\cdot d_2 - \sum \nu_{ij}, \]
where the sum runs over linearly independent rows of ${\mathcal{B}}$ lying in opposite
open quadrants of~$\Z^2$. But then, since
\[\deg(\ind_w(H_{\mathcal{B}}(c))) = \rank(\ini_{(-w,w)}(H_{\mathcal{B}}(c))
\leq \rank(H_{\mathcal{B}}(c))\]
we conclude that
\[ \rank(H_{{\mathcal{B}}}(c)) = \rank(\Horn({{\mathcal{B}}},c) ) \geq
d_1 \cdot d_2 - \sum \nu_{ij}.\]
The reverse inequality follows from Corollary~\ref{coro:horn-rank-bound}.

\end{proof}

The same method exactly that proved Theorem~\ref{thm:rank-formula} will
compute the rank of the hypergeometric system arising from a
lattice (actually, this proof is easier, since 
$\# T(\ini_w(I_{\mathcal{B}})) = \deg(I_{\mathcal{B}}) = g \cdot \vol(A)$
is easier to compute than $\# T(\ini_w(I))$).
Notice that here we do not need to require that $m=2$, since
we know what the solutions of these systems look like without restriction
on the codimension of~$I_{\mathcal{B}}$.

\begin{theorem}

For generic~$c$,
\[ \rank ( I_{\mathcal{B}} +\langle A \cdot \theta - A\cdot c \rangle ) =\#
T(\ini_w(I_{\mathcal{B}})) = \deg(I_{\mathcal{B}}) = g \cdot \vol(A) .\]
\end{theorem}


\section{Explicit construction of fully supported hypergeometric functions}
\label{sec:combin}

We already know how to explicitly write down Puiseux polynomial solutions of
a bivariate Horn system with generic parameters. 
This is done by taking pairs of rows of the matrix~${\mathcal{B}}$
that are linearly independent and lie in opposite open quadrants of~$\Z^2$,
obtaining a cone from these vectors, and joining together lattice points in the cone
using horizontal and vertical moves to obtain the finite supports of Puiseux
polynomial solutions. We have not described the coefficients appearing in these Puiseux
polynomials, although they are easily computed on a case by case basis.

The goal of this section is to be even more explicitly describe the fully
supported solutions of~${H_{\mathcal{B}}(c)}$, and thus of $\Horn({\mathcal{B}},c)$.
In particular, we will show that the fully supported solutions of
$\Horn({\mathcal{B}},c)$ are hypergeometric in the following classical sense.

\begin{definition}
\label{def:hyp}
A formal power series $\sum_{(s,t)\in \Z^2} \lambda(s,t) y_1^s y_2^t$ is 
{\em hypergeometric} if there exist rational functions~$R_1$ and~$R_2$ such that:
\[ \lambda(s+1,t) = R_1(s,t)\lambda(s,t) \quad \mbox{and} \quad
\lambda(s,t+1) = R_2(s,t)\lambda(s,t). \]
In this paper we restrict our attention to the case when the numerator and 
the denominator
of the rational functions~$R_1,R_2$ are products of affine linear functions with 
integer coefficients by~$s,t$ and arbitrary constant terms. 
\end{definition}

A formal power series such as in Definition~\ref{def:hyp} satisfies a Horn system 
of differential equations. We will now show that the other
fully supported solutions of this system are spanned by monomial multiples of series of this form.
We know that the fully supported solutions of~$H_{\mathcal{B}}(c)$ are simply
the solutions of the lattice hypergeometric system
$I_{\mathcal{B}} + \langle A\cdot \theta - A\cdot c \rangle$. The following
result is proved using the methods from~\cite[Section 3.4]{SST}. We start
by setting up some notation. Recall that~$L_{\mathcal{B}}$ is the
lattice in~$\Z^n$ spanned by the columns of~${\mathcal{B}}$.

Given $v \in \C^n$ we let
\[ N_v = \{ u \in L_{{\mathcal{B}}}:
v_i \in \Z_{<0} \Leftrightarrow (u+v)_i \in \Z_{<0} \quad \mbox{and} \quad
v_i \in \Z_{\geq 0} \Leftrightarrow (u+v)_i \in \Z_{\geq 0}\},\]
and define a formal power series
\begin{equation}
\label{eqn:phiv}
\phi_v := x^v \sum_{u \in N_v} \frac{[v]_{u_-}}{[v+u]_{u_+}} x^u ,
\end{equation}
where
\[ [v]_{u_-} = \prod_{i:u_i<0} \prod_{j=1}^{-u_i} (v_i-j+1) \quad
\mbox{and} \quad [v+u]_- = \prod_{i:u_i>0} \prod_{j=1}^{u_i} (v_i+j). \]

\begin{theorem}
\label{thm:form-of-sols}
Let~$c$ be generic and~$w$ a generic weight vector.
Denote by $v^{(1)},\dots,v^{(g \cdot \vol(A))}$ be the zeros of the indicial ideal
$\ind_w(I_{\mathcal{B}}+\langle A\cdot \theta - A\cdot c\rangle)$. Then
the formal power series $\{\phi_{v^{(i)}} : i=1,\dots,g\cdot \vol(A)\}$
are linearly independent holomorphic solutions of $I_{\mathcal{B}}+
\langle A\cdot \theta - A\cdot c \rangle$.
\end{theorem}

\begin{proof}
For sufficiently generic~$c$, the vectors~$v^{(i)}$ have no negative 
integer coordinates.
Now use the arguments from~\cite[Theorem 3.4.2]{SST}. In particular,
the support of each of these series is contained in a strongly convex cone.
\end{proof}

We now have an explicit description of a basis of the solution space of
$\Horn({{\mathcal{B}}},c)$ (and~$H_{{\mathcal{B}}}(c)$).

\begin{theorem}

\label{thm:sol-basis}
If~$c$ is generic, the fully supported
series obtained by applying the isomorphism
from Corollary~\ref{coro:vector-space-iso} to the fully supported series
constructed in Theorem~\ref{thm:form-of-sols}
and the Puiseux polynomials
constructed in Theorem~\ref{thm:puiseux}
form a basis for the solution space of $\Horn({{\mathcal{B}}},c)$.
\end{theorem}

\begin{proof}

Theorem~\ref{thm:form-of-sols} and  Corollary
\ref{coro:vector-space-iso} give us $g \cdot \vol(A) + \sum \nu_{ij}$
linearly independent solutions of $H_{\mathcal{B}}(c)$ (here the sum
runs over linearly independent rows of~${\mathcal{B}}$). By Theorem
\ref{thm:rank-formula}, these must span the solution space of~$H_{\mathcal{B}}(c)$.
\end{proof}

Notice that applying the change of variables from Corollary~\ref{coro:vector-space-iso}
to the functions~$\phi_{v^{(i)}}$ from Theorem~\ref{thm:form-of-sols}
is particularly easy.

\begin{corollary}
For~$c$ generic and~$v^{(i)}$ as in Theorem~\ref{thm:form-of-sols},
let~$\alpha^{(i)}$ be the unique vector that satisfies $v^{(i)}-c={\mathcal{B}} \cdot
\alpha^{(i)}$.
Then the space
of fully supported solutions of $\Horn({\mathcal{B}},c)$ is spanned by the functions
\[y_1^{\alpha^{(i)}_1} y_2^{\alpha^{(i)}_2}
\sum_{{\mathcal{B}}\cdot z \in N_{v^{(i)}}}
\frac{[v^{(i)}]_{({\mathcal{B}}\cdot z)_-}}{
[v^{(i)}+{\mathcal{B}}\cdot z]_{({\mathcal{B}}\cdot z)_+}}
y_1^{z_1} y_2^{z_2} .\]
In particular, all the fully supported solutions of $\Horn(\mathcal{B},c)$
are spanned by  monomial multiples of hypergeometric series in the sense of 
Definition~\ref{def:hyp}.
\end{corollary}


\section{Holonomicity of $\Horn({{\mathcal{B}}},c)$}

\label{sec:holonomic}

Throughout this section we assume that $m=2$.
Since we do not have a $D$-module isomorphism between $H_{{\mathcal{B}}}(c)$ and
$\Horn({\mathcal{B}},c)$, the holonomicity of $H_{{\mathcal{B}}}(c)$ does not
directly prove that $\Horn({\mathcal{B}},c)$ is holonomic.
In this section we prove that the bivariate hypergeometric system
$\Horn({\mathcal{B}},c)$
is holonomic.

Recall that a system of differential equations is said to be {\it holonomic}
if the dimension of its characteristic variety is the same as
the dimension of the variable space.

We recall that we are dealing with the system of equations defined
by the hypergeometric operators
\begin{equation}
\begin{array}{l}
H_{1} = \bm{Q}_{1}(\theta) - y_{1}\bm{P}_{1}(\theta),   \\
H_{2} = \bm{Q}_{2}(\theta) - y_{2}\bm{P}_{2}(\theta).
\end{array}
\label{horn}
\end{equation}
By the definition of the Horn system (see Section~\ref{sec:formula})
the bivariate polynomials~$\bm{P}_{i},\bm{Q}_{i}$ satisfy the compatibility condition
\begin{equation}
R_{1}(s + e_2)R_{2}(s) = R_{2}(s + e_1)R_{1}(s),
\label{compatibility}
\end{equation}
where $R_{i}(s) = \bm{P}_{i}(s)/\bm{Q}_{i}(s + e_i)$ and $\{ e_1,e_2 \}$ is
the standard basis of~$\Z^2.$

\begin{theorem}

A bivariate Horn system with generic parameters is holonomic.
\label{hornholonom}

\end{theorem}

To prove this theorem we need some intermediate results and notation.
Let~$(H_{1},H_{2}) \subset D_2$ denote the ideal generated by
the hypergeometric operators defining the Horn system.
By~$\sigma(P)$ we denote the principal symbol of a differential
operator~$P$. This is an element of the polynomial ring $\C[y_1,y_2,z_1,z_2]$.
The only case when a bivariate Horn system is not holonomic
is when the principal symbols of all the operators in~$(H_{1},H_{2})$
have a nontrivial greatest common divisor (for otherwise we have two
independent algebraic equations and hence the dimension of the
characteristic variety of the Horn system is~2).
Thus to prove holonomicity of~(\ref{horn}) it suffices to construct
a family of operators in~$(H_{1},H_{2})$ such that the greatest common
divisor of their principal symbols is~$1.$

By the construction of the operators in the Horn system
(see Section~\ref{sec:formula}) the greatest common
divisor of the principal symbols of~$H_{1}$ and~$H_{2}$ is given by
a product of powers of linear forms $a y_{1}z_{1} + b y_{2}z_{2},$
where~$a,b\in\Z.$
Thus to prove Theorem~\ref{hornholonom} it suffices to show that for
any~$a,b\in\Z$ such that $a y_{1}z_{1} + b y_{2}z_{2}$ divides
${\rm gcd}(\sigma(H_{1}),\sigma(H_{2}))$ there exists an
operator~$T_{a,b}\in (H_{1},H_{2})$ whose principal symbol is not
divisible by $a y_{1}z_{1} + b y_{2}z_{2}.$

\begin{remark}

\rm
For generic parameters the compatibility
condition~(\ref{compatibility}) is equivalent to the relations
\begin{equation}
[ y_{1}\bm{P}_{1}(\theta), y_{2}\bm{P}_{2}(\theta) ] = 0, \quad
(E_{2}\bm{Q}_{2})(\theta) (E_{1}E_{2}\bm{Q}_{1})(\theta) =
(E_{1}\bm{Q}_{1})(\theta) (E_{1}E_{2}\bm{Q}_{2})(\theta),
\label{diffcomp}
\end{equation}
where~$[\, , ]$ denotes the commutator of two operators,
$(E_{i}^{\lambda} P)(s) = P(s + \lambda e_{i})$
and~$E_{i} = E_{i}^{1}.$
Indeed, the equalities~(\ref{diffcomp}) mean that
the numerators (respectively the denominators) of the rational
functions in~(\ref{compatibility}) are equal. The generic parameters
assumption implies that no cancellations can occur and hence this is
indeed the case.
\end{remark}

\begin{lemma}

For any $\alpha,\beta,\gamma,\delta \in \C$
and ${P}_1(\theta),{P}_2(\theta),{Q}_1(\theta),{Q}_2(\theta)$
satisfying the relations
\begin{equation}
[ y_{1}{P}_{1}(\theta), y_{2}{P}_{2}(\theta) ] = 0, \quad
(E_{2}{Q}_{2})(\theta) (E_{1}E_{2}{Q}_{1})(\theta) =
(E_{1}{Q}_{1})(\theta) (E_{1}E_{2}{Q}_{2})(\theta),
\label{diffcomp2}
\end{equation}
it holds that:
\begin{align}
\label{therelation}
(\alpha (E_{2}^{-1} {Q}_{1})(\theta) - \beta y_{1}{P}_{1}(\theta))
&(\gamma {Q}_{2}(\theta) - \delta y_{2}{P}_{2}(\theta)) - \\ \nonumber
&(\alpha (E_{1}^{-1}{Q}_{2})(\theta) - \beta y_{2}{P}_{2}(\theta))
(\gamma {Q}_{1}(\theta) - \delta y_{1}{P}_{1}(\theta)) =
\left|
\begin{array}{cc}
\alpha &  \beta   \\
\gamma &  \delta  \\
\end{array}
\right| \Psi,
\end{align}
where
$\Psi= y_{1}{Q}_{2}(\theta){P}_{1}(\theta) - y_{2}{Q}_{1}(\theta){P}_{2}(\theta).$
\label{psilemma}

\end{lemma}

The proof of Lemma~\ref{psilemma} is a direct computation which uses
the compatibility conditions~(\ref{diffcomp2})
and the Weyl algebra identity
$(E_{i}^{-1}{Q}_{j})(\theta) y_{i} = y_{i} {Q}_{j}(\theta).$

Let us now consider a special case to which we will later reduce the case of
an arbitrary bivariate Horn system with generic parameters. Namely, let us
find a holonomicity condition for the system defined by the operators:
\begin{equation}
\begin{array}{l}
U_{1} = f(t) {Q}_{1}(\theta) - y_{1} g(t) {P}_{1}(\theta),  \\
U_{2} = f(t) {Q}_{2}(\theta) - y_{2} g(t) {P}_{2}(\theta),
\end{array}
\label{specialhorn}
\end{equation}
where~$f,g$ are arbitrary non zero univariate polynomials, $t=\theta_{1} + \theta_{2}$
and~${P}_{i},{Q}_{i}$ are arbitrary bivariate polynomials such
that
${\rm deg}\, f + {\rm deg}\, {Q}_{i} = {\rm deg}\, g + {\rm deg}\, {P}_{i}$ and
that~${P}_{i},{Q}_{i}$ satisfy~(\ref{diffcomp2}). Note that these relations are satisfied
if $g(t) P_i, g(t) Q_i$ satisfy the equivalent relations.
We assume also that~$t$ is not
present in~${P}_{i}(\theta), {Q}_{i}(\theta),$ i.e., that none of the principal
symbols of these operators vanish along the hypersurface
$y_1 z_1 + y_2 z_2 = 0.$

Our goal is to ``eliminate~$t$'' from~(\ref{specialhorn}), i.e.,
to construct an operator in the ideal~$(U_{1},U_{2})$ whose principal
symbol is not divisible by $\sigma(t)=y_1 z_1 + y_2 z_2.$
We do it as follows.

\begin{lemma} \label{psilemma2}
Let~$\Psi$ be as in Lemma~\ref{psilemma}. Then
$ R(f(t),g(t)) \Psi  \in (U_{1},U_{2})$, where~$R(f(t),g(t))$ is the
resultant of~$f,g.$
\label{resintheideal}
\end{lemma}

\begin{proof}

Let us write the polynomials~$f,g$ in the form
$f(t) = \sum_{i=0}^{d} f_{i} t^{i},$ $g(t) = \sum_{i=0}^{d} g_{i} t^{i}.$
Notice that~$f,g$ do not have to be of the same degree since some
of~$f_{i},d_{i}$ may be zero. Using~(\ref{therelation}),
and the fact that the subring of the Weyl algebra generated
by~$\theta_1$ and~$\theta_2$ is commutative, we conclude
that for any~$j=0,\ldots,d$
\begin{align*}
\sum_{i=0}^{d}
\left|
\begin{array}{cc}
f_{j} & g_{j} \\
f_{i} & g_{i} \\
\end{array}
\right|
\Psi t^{i} & = \sum_{i=0}^d \Theta_{1j} (f_it^i {Q}_{2}(\theta) - y_2 g_i t^i {P}_{2}(\theta)) -
\sum_{i=0}^d \Theta_{2j} (f_it^i {Q}_{1}(\theta)-y_1g_it^i{P}_{1}(\theta)) \\
& = \Theta_{1j} U_2 + \Theta_{2j} U_1 \in (U_1,U_2),
\end{align*}
where $\Theta_{1j} = f_j (E_{2}^{-1}{Q}_{1})(\theta) - g_j y_{1}{P}_{1}(\theta)$ and
$\Theta_{2j} = f_j (E_{1}^{-1}{Q}_{2})(\theta) - g_j y_{2}{P}_{2}(\theta)$.

Now clearly,
\[\sum_{i=0}^{d}
\left|
\begin{array}{cc}
f_{j} & g_{j} \\
f_{i} & g_{i} \\
\end{array}
\right|
\Psi t^{i} =
\Psi
\left|
\begin{array}{cc}
f_j  &  g_j   \\
f(t) &  g(t)  \\
\end{array}
\right|,
\]
so that:
\[ \Psi
\left|
\begin{array}{cc}
f_j  &  g_j   \\
f(t) &  g(t)  \\
\end{array}
\right| \in (U_1,U_2).\]
In the trivial case when the polynomials~$f$ and~$g$ are proportional we
have~$R(f(t),g(t))=0$ and the conclusion of the lemma is obviously true.
If~$f$ is not proportional to~$g$ then the rank of the
$2 \times (m + 1)$-matrix
$
\left(
\begin{array}{ccc}

f_{0} \ldots f_{m} \\
g_{0} \ldots g_{m}

\end{array}
\right)
$
equals~$2$ and hence $\Psi f(t), \Psi g(t) \in (U_{1},U_{2}).$
Since $(t-1)\Psi = \Psi t,$ it follows that
$\Psi h(t) \in (U_{1},U_{2}),$ for any~$h(t) \in (f(t),g(t)),$
where~$(f(t),g(t))$ denotes the ideal in the ring of (commuting)
univariate polynomials generated by~$f,g.$ It is known that the resultant
of two polynomials lies in the ideal generated by these polynomials and
hence $R(f(t),g(t)) \Psi \in (U_{1},U_{2}).$
The proof is complete.
\end{proof}

\begin{corollary}

Suppose that ${\rm gcd} (\sigma(U_{1}), \sigma(U_{2}))$ is a power
of $x_{1}z_{1} + x_{2}z_{2}.$
Then the hypergeometric system~(\ref{specialhorn}) is holonomic if and only if
$R(f(t),g(t)) \neq 0.$
\label{specialholonom}

\end{corollary}
\begin{proof}

Suppose that~$R(f(t),g(t)) = 0$ and let~$\zeta\in\C$ be a common root
of the polynomials~$f,g.$ Since for any smooth univariate function~$h$
the product~$y_{2}^{\zeta} h(y_{1}/y_{2})$ is annihilated by the
operator $t-\zeta = \theta_{1} + \theta_{2} - \zeta,$ it follows that
the space of analytic solutions to~(\ref{specialhorn}) has infinite
dimension.
It is known that a holonomic system can only have finitely many
linearly independent solutions and hence~(\ref{specialhorn}) is not
holonomic in this case.

On the other hand, if~$R(f(t),g(t))\neq 0,$ then by
Lemma~\ref{resintheideal} the operator~$\Psi$
lies in the ideal~$(U_{1},U_{2}).$ By the assumption of the corollary
the principal symbols of~$U_{1},U_{2}$ and~$\Psi$ are relatively prime
and hence the system~(\ref{specialhorn}) is holonomic.
\end{proof}

\begin{example}
Consider the system quoted in the introduction, given by the two
hypergeometric operators
$$
H_1 \,  =  \,
x (\theta_x + \theta_y +a)(\theta_x+b)-\theta_x(\theta_x+\theta_y+c-1) , \\
H_2 \,  =  \,  y (\theta_x + \theta_y +a)(\theta_y+b')-\theta_y(\theta_x+\theta_y+c-1)
$$
for  Appell's function~$F_1$.
The operator~$\Psi$ in Lemma~\ref{psilemma} equals in this case
$$\Psi =  (x \ y)  \  \Psi'\ , \quad  {\rm where}  \quad \Psi' =  (x-y) \partial_x \partial_y
-b' \partial_x + b \partial_y.$$
When $a-c+1 \not=0,$ we deduce from Lemma~\ref{psilemma2} that $(x\ y) \ \Psi'$ lies
in the $D$-ideal  $\langle H_1, H_2 \rangle.$ In particular,
all holomorphic solutions~$\varphi$ of the Appell system will also satisfy
$\Psi'(\varphi) =0.$ We point out that some authors add this third equation to the
system (cf. for instance~\cite[Page 48]{SST}). In fact,
having this operator, the holonomicity of the system follows immediately.
\end{example}

We are now in a position to complete the proof of
Theorem~\ref{hornholonom}.

\vskip0.1cm

\begin{proof}[Proof of Theorem~\ref{hornholonom}]
Suppose that the polynomial ${\rm gcd} (\sigma(H_{1}),\sigma(H_{2}))$ vanishes along
the hypersurface $a y_{1}z_{1} + b y_{2}z_{2} = 0.$ We aim to
construct an operator in the ideal~$(H_{1},H_{2})$ whose principal
symbol is not divisible by $a y_{1}z_{1} + b y_{2}z_{2}.$
The change of variables $\xi_{1} = y_{1}^{1/a},$ $\xi_{2} = y_{2}^{1/b}$
transforms the operator $a \theta_{y_{1}} + b \theta_{y_{2}}$ into
the operator $\theta_{\xi_{1}} + \theta_{\xi_{2}}$ and the
system~(\ref{horn}) into the system generated by the operators
\begin{equation}
\begin{array}{l}
\hat{{Q}}_{1}(\theta_{\xi_{1}},\theta_{\xi_{2}}) -
\xi_{1}^{a} \, \hat{{P}}_{1}(\theta_{\xi_{1}},\theta_{\xi_{2}}),   \\
\hat{{Q}}_{2}(\theta_{\xi_{1}},\theta_{\xi_{2}}) -
\xi_{2}^{b} \, \hat{{P}}_{2}(\theta_{\xi_{1}},\theta_{\xi_{2}}),
\end{array}
\end{equation}
where~$\hat{{P}}_{i}(u,v)={{P}}_{i}(u/a,v/b),$
$\hat{{Q}}_{i}(u,v)={Q}_{i}(u/a,v/b).$

Let us introduce operators $\lambda_{ia}^{k},\mu_{ia}^{k}$ acting on a
bivariate polynomial~$P$ as follows:
\begin{equation}
\lambda_{ia}^{k} (P) = \prod_{j=1}^{k}   (E_{i}^{-ja} P), \qquad
\mu_{ia}^{k}     (P) = \prod_{j=0}^{k-1} (E_{i}^{ja}  P).
\label{adjustop}
\end{equation}
(Notice that the upper index here is {\bf not} a power.)
The next Weyl algebra identities follow directly from the definition
of~$\lambda_{ia}^{k},\mu_{ia}^{k}$ (the arguments of all of the
involved polynomials being~$\theta_{\xi_{1}}, \theta_{\xi_{2}}$):
\begin{equation}
\begin{array}{rcl}
\lambda_{1a}^{k} \, (\hat{{Q}}_{1})  \xi_{1}^{a} &=&
\xi_{1}^{a} \, \hat{{Q}}_{1} \,
\lambda_{1a}^{k-1} \, (\hat{{Q}}_{1}), \\
\lambda_{2b}^{k} \, (\hat{{Q}}_{2})  \xi_{2}^{b} &=&
\xi_{2}^{b} \, \hat{{Q}}_{2} \,
\lambda_{2b}^{k-1} \, (\hat{{Q}}_{2}), \\
\mu_{1a}^{k} \, (\hat{{P}}_{1})  \xi_{1}^{a} \,
\hat{{P}}_{1} &=&
\xi_{1}^{a} \, \mu_{1a}^{k+1} \, (\hat{{P}}_{1}), \\
\mu_{2b}^{k} \, (\hat{{P}}_{2}) \xi_{2}^{b} \,
\hat{{P}}_{2} &=&
\xi_{2}^{b} \, \mu_{2b}^{k+1} \, (\hat{{P}}_{2}).
\end{array}
\label{adjustopiden}
\end{equation}

Using~(\ref{adjustopiden}) we arrive at the equalities
\begin{align}
\label{transformedhorn1}
\left(
\sum_{\nu=0}^{b-1} \xi_{1}^{\nu a} \lambda_{1a}^{b-1-\nu}
(\hat{{Q}}_{1})(\theta_{\xi}) \mu_{1a}^{\nu}(\hat{{P}}_{1})(\theta_{\xi})
\right)
(\hat{{Q}}_{1}(\theta_{\xi}) - & \xi_{1}^{a} \hat{{P}}_{1}(\theta_{\xi})) = \\ \nonumber
& \hat{{Q}}_{1}(\theta_{\xi}) \lambda_{1a}^{b-1}
(\hat{{Q}}_{1})(\theta_{\xi})  -
\xi_{1}^{ab} \mu_{1a}^{b} (\hat{{P}}_{1})(\theta_{\xi}),
\end{align}

\begin{align}
\label{transformedhorn2}
\left(
\sum_{\nu=0}^{a-1} \xi_{2}^{\nu b} \lambda_{2b}^{a-1-\nu}
(\hat{{Q}}_{2})(\theta_{\xi}) \mu_{2b}^{\nu}(\hat{{P}}_{2})(\theta_{\xi})
\right)
(\hat{{Q}}_{2}(\theta_{\xi}) - & \xi_{2}^{b} \hat{{P}}_{2}(\theta_{\xi})) = \\ \nonumber
&\hat{{Q}}_{2}(\theta_{\xi}) \lambda_{2b}^{a-1}
(\hat{{Q}}_{2})(\theta_{\xi})  -
\xi_{2}^{ab} \mu_{2b}^{a} (\hat{{P}}_{2})(\theta_{\xi}).
\end{align}

The differential
operators~(\ref{transformedhorn1}) and (\ref{transformedhorn2}) are
Horn-type hypergeometric operators in the
variables~$\eta_{1}= \xi_{1}^{ab}$ and~$\eta_{2}=\xi_{2}^{ab}.$
Let us write these operators in the form
$$
\begin{array}{l}
\tilde{U}_{1} = f(\tau) \tilde{{Q}}_{1}(\theta_{\eta}) -
\eta_{1} g(\tau) \tilde{{P}}_{1}(\theta_{\eta}), \\
\tilde{U}_{2} = f(\tau) \tilde{{Q}}_{2}(\theta_{\eta}) -
\eta_{2} g(\tau) \tilde{{P}}_{2}(\theta_{\eta}), \\
\end{array}
$$
where~$f,g$ are univariate polynomials,
$\tau=\theta_{\eta_{1}} + \theta_{\eta_{2}}$ and none of the
principal symbols of the operators
$\tilde{{P}}_{i}(\theta_{\eta}), \tilde{{Q}}_{i}(\theta_{\eta})$
vanish along the hypersurface $\eta_{1} z_{1} + \eta_{2} z_{2}=0.$
The existence of such polynomials~$f,g$
follows from the compatibility condition which is satisfied
by~(\ref{transformedhorn1}),(\ref{transformedhorn2}).

By Lemma~\ref{resintheideal} the operator
$\tilde{\Psi} = \eta_{1}\tilde{{Q}}_{2}(\theta_{\eta}) \tilde{{P}}_{1}(\theta_{\eta})-
\eta_{2}\tilde{{Q}}_{1}(\theta_{\eta}) \tilde{{P}}_{2}(\theta_{\eta})$
lies in the ideal~$(\tilde{U}_{1},\tilde{U}_{2})$ as long as the
parameters of the original Horn system~(\ref{horn}) are generic.
Notice that by construction the principal symbol of~$\tilde{\Psi}$
does not vanish along the hypersurface $\eta_{1} z_{1} + \eta_{2} z_{2}=0.$
Going back to the variables~$y_{1},y_{2},$ we conclude that there exists
an operator in~$(H_{1},H_{2})$ whose principal symbol is not divisible
by $a y_{1} z_{1} + b y_{2} z_{2}.$
This completes the proof of
Theorem~\ref{hornholonom}.
\end{proof}


\section{The Cohen-Macaulay property as a tool to compute rank, and
  further research directions}

Since the lattice basis ideal~$I$ is a complete intersection and therefore
Cohen-Macaulay, it is natural to try to apply the methods that proved that
the holonomic rank $H_A(A\cdot c)$ is always~$\vol(A)
=\deg(I_A)$ when the underlying
toric ideal~$I_A$ is Cohen-Macaulay.

The first evidence that these methods will not work is that the generic rank
of the Horn system $H_{{\mathcal{B}}}(c)$ is {\em not}
$\deg(I)= d_1 \cdot d_2$, unless we make the assumption that~${{\mathcal{B}}}$
has no linearly dependent rows in opposite open quadrants of~$\Z^2$.

If we follow the arguments that proved~\cite[Lemma 4.3.7]{SST}, which is the main
ingredient needed to prove that, when~$I_A$ is Cohen-Macaulay,
$\rank(H_A(A\cdot c)) = \vol(A)$ for all~$c$, we see that the crucial point is whether
the $n-m$ polynomials
\begin{equation}
\label{eqn:reg-seq}
\sum_{j=1}^n a_{ij} x_j z_j \in \C[x_1,\dots ,x_n,z_1,\dots ,z_n],\; i=1,
\dots ,n-m,
\end{equation}
form a regular sequence in $\C(x_1,\dots ,x_n)[z_1, \dots ,z_n]/I$, where
here we think of~$I$ as an ideal in the variables $z_1, \dots,z_n$.
But if~${{\mathcal{B}}}$ has linearly dependent rows in opposite open quadrants,
the ring
\[ (\C(x_1,\dots ,x_n)[z_1,\dots ,z_n])/(I+ \langle
\sum_{j=1}^n a_{ij} x_j z_j \in \C[x_1,\dots ,x_n,z_1,\dots ,z_n],\; i=1,
\dots ,n-m \rangle)\]
is not artinian!

\begin{lemma}
\label{lemma:hsop}
Let $m=2$.
If $\langle A \cdot xz \rangle$ is ideal generated by the polynomials (\ref{eqn:reg-seq}),
then the ideal $I+\langle A\cdot xz\rangle$ is artinian
in $\C(x_1,\dots ,x_n)[z_1,\dots ,z_n]/I$, if and only if~${{\mathcal{B}}}$ 
has no linearly dependent rows in opposite open quadrants of~$\Z^2$.
\end{lemma}

\begin{proof}

We need to investigate the intersection of the zero locus of $\langle A
\cdot xz \rangle$ over $\C(x)$ with the zero locus of~$I$ over~$\C(x)$. 
Specifically, we want to show that this intersection is a
finite set if and only if~${{\mathcal{B}}}$ contains no linearly
dependent rows in opposite open quadrants of~$\Z^2$. We can perform
this intersection irreducible component by irreducible component of~$I$, 
recalling the primary decomposition of~$I$ from 
Proposition~\ref{propo:lattice-basis-description}.

The toric irreducible components of~$I$ we can deal with all at the
same time: we know that $\C(x)[z]/(I_{{\mathcal{B}}}+\langle A\cdot xz
\rangle)$ is zero-dimensional. That just leaves the primary components
of~$I$ corresponding to associated primes $\langle z_i, z_j \rangle$,
where~$b_i$ and~$b_j$ lie in the interior of open quadrants of~$\Z^2$. 
But now it is clear that such a component will meet the zero
locus of $\langle A\cdot xz \rangle$ in an infinite set if and only if~$b_i$ 
and~$b_j$ are linearly dependent.
\end{proof}

As a consequence of Lemma~\ref{lemma:hsop} and the arguments in~\cite[Section 4.3]{SST},
we have one case when the fact that~$I$ is a complete intersection will imply that
the rank of~$H_{{\mathcal{B}}}(c)$ does not depend on~$c$.

\begin{theorem}

If~${{\mathcal{B}}}$ has no linearly dependent rows in opposite quadrants of~$\Z^2$ then
\[ \rank(H_{{\mathcal{B}}}(c)) = d_1\cdot d_2 \quad \mbox{for all}\; c \in \C^n.\]
\end{theorem}

Notice that this result holds even when the rows of~${\mathcal{B}}$ do
not add up to zero. 

Remark that the case in which no pair of (linearly dependent or not) rows lie
in the interior of
opposite quadrants corresponds precisely to the case in which
the lattice ideal~$I_{\mathcal{B}}$ is a complete
intersection. This agrees with the characterization in~\cite{fsh}.

There is another situation when we can apply the arguments from
\cite[Section 4.3]{SST} to prove that a certain holonomic rank does not depend on~$c$. 
Let~$J$ be the ideal in $\C[\partial_1,\dots,\partial_n]$ obtained by
saturating from~$I$ the components~$I_{ij}$ corresponding to
linearly dependent rows of~${{\mathcal{B}}}$. Then
\[ \deg(J)= d_1\cdot d_2 -\sum \nu_{ij} ,\]
where the sum runs over the linearly dependent rows of~${{\mathcal{B}}}$ that
lie in opposite open quadrants of~$\Z^2$. As before, the methods 
in~\cite[Section 4.3]{SST} prove the following result.

\begin{lemma}

If~$J$ is Cohen-Macaulay,
\[ \rank (J+\langle A \cdot \theta - A\cdot c \rangle) = \deg(J) . \]
\end{lemma}

The previous lemma and our rank formula for Horn systems have the
following consequence.

\begin{corollary}

\label{coro:too-strong-assumptions}
If~$J$ is Cohen-Macaulay and~$c$ is generic, the solution spaces 
of~$H_{\mathcal{B}}(c)$ and $J+\langle A\cdot \theta - A\cdot c \rangle$ coincide.
\end{corollary}

We believe that Corollary~\ref{coro:too-strong-assumptions} holds even
when~$J$ is not Cohen-Macaulay. It would be desirable to obtain an
independent proof of this, since in that case we would have a proof of
our rank formula in the case that~$J$ is Cohen-Macaulay that does not
rely on a precise description of the solution space. 

\bigskip

The natural question at this point is whether we can extend arguments in 
Section~\ref{sec:bound} to give an algebraic
formula for the rank of a Horn system for any~$m$. However, in order to
use those methods, several ingredients are missing. First, we need to assume that
the lattice basis ideal~$I$ is a complete intersection, since this is not necessarily
true if~$m>2$.
Moreover, it is not true in general that given a toric ideal~$I_A$, one can find
a lattice basis ideal contained in~$I_A$ that is a complete intersection~\cite{EC}.
Moreover,
our techniques for finding the form of the solutions of~$H_{\mathcal{B}}(c)$
for~$m=2$ do not directly generalize to higher~$m.$
In any case, in order to obtain an explicit rank formula in the case
that $m>2$,  combinatorial expressions for the multiplicities
of the minimal primes of any lattice basis ideal are needed.

\bigskip
\noindent {\bf Acknowledgments:}
Work on this article began during T.~Sadykov's visit to
the University of Buenos Aires in September~2001 and a
subsequent meeting of the three authors in Stockholm in January~2002. 
We are very grateful to Mikael Passare,
who made those visits possible. We thank Jan-Erik
Bj\"{o}rk for all his help and patience and Teresa Monteiro Fernandes for
insights on regular holonomic D-modules.
We would also like to thank Bernd Sturmfels  and
Michael Singer for inspiring conversations. We are also grateful
to Ezra Miller, who suggested the proof of Lemma~\ref{lemma:multiplicities}.
Part of this work was completed while the first two authors
were members at~MSRI, which we also thank for its support and
wonderful working atmosphere.

\bibliographystyle{plain}

\def\cprime{$'$} \def\cprime{$'$}

\end{document}